\documentclass[11pt]{article}
\usepackage{amsmath}
\usepackage{amsfonts}
\usepackage{amssymb}
\usepackage{amsthm}
\usepackage{graphicx}

\theoremstyle{plain}
\newtheorem{theorem}{Theorem}[section]
\newtheorem{proposition}[theorem]{Proposition}
\newtheorem{lemma}[theorem]{Lemma}
\newtheorem{corollary}[theorem]{Corollary}
\theoremstyle{remark}
\newtheorem{remark}[theorem]{Remark}
\numberwithin{equation}{section}
\begin{document}

\title{Neumann Series of Bessel Functions\\
[-.2ex] for a Quadratic Sturm--Liouville Pencil}
\author{Vladislav V. Kravchenko$^{1}$, Sergii M. Torba$^{1}$ \\
{\small $^{1}$ Cinvestav, Unidad Quer\'{e}taro, Departamento de Matem\'{a}ticas}\\
{\small Libramiento Norponiente \#2000, Fracc. Real de Juriquilla,}\\
{\small Quer\'{e}taro, Qro., 76230 M\'{e}xico.}\\
{\small e-mail: vkravchenko@math.cinvestav.edu.mx}\\
{\small e-mail: storba@math.cinvestav.edu.mx, sergii.torba@cinvestav.mx}}
\date{}
\maketitle

\begin{abstract}
We construct Neumann series of Bessel functions representations for 
solutions and their first derivatives of the equation $-y^{\prime\prime}(x)+(q(x)+2\rho p(x))y(x)=\rho^2y(x)$. The series arise from
Legendre expansions of two transmutation operator kernels. The obtained series converge uniformly on horizontal strips in the complex plane of the parameter $\rho$. We develop a recursive integration procedure for computing the coefficients of these representations. The knowledge of the first coefficients allows one to recover the potentials $p$ and $q$ which is useful for solving inverse problems. For separated boundary conditions, we
estimate the characteristic functions whose zeros determine the eigenvalues.
Two numerical examples illustrate the practical applicability of the obtained series representations to the solution of spectral problems.
\end{abstract}



\section{Introduction}

We consider the quadratic Sturm--Liouville pencil on a finite interval
\begin{equation}
-y^{\prime \prime }(x)+\left( q(x)+2\rho p(x)\right) y(x)=\rho
^{2}y(x),\qquad 0<x<b.  \label{eq:pencil}
\end{equation}%
Here $\rho \in \mathbb{C}$ is a spectral parameter, and the complex-valued
coefficients satisfy $q\in L^{2}(0,b)$ and $p\in W^{1,2}(0,b)$.

The aim of the present work is to derive new series representations for the
solutions of (\ref{eq:pencil}) and for their derivatives in the form of
Neumann series of Bessel functions (NSBF). This extends the result of \cite%
{KNT2017} to quadratic Sturm--Liouville pencils and provides a basis for
methods addressing direct and inverse problems related to (\ref{eq:pencil}).

Quadratic pencils arise in energy dependent Schr\"{o}dinger equations and in
inverse spectral problems for diffusion operators \cite%
{GasymovGuseinov1981,Guseinov1986,Guseinov2014,GuseinovNabiev2007}. Related
models describe vibrations of strings with viscous drag and Klein--Gordon
evolution \cite{Yamamoto1990,LangerNajmanTretter2008,HrynivPronska2012}. The
inverse scattering problem for an energy dependent potential was studied in
\cite{JaulentJean1976I}. Direct and inverse spectral problems for quadratic pencils are studied in \cite{ButerinYurko2012, LyuWei2025, Pronska2013, Pronska2013TwoSpectra, Pronska2013Spectral, Yurko2000Pencils, Yurko2002, Yurko2012Quasiperiodic}. Inverse problems with complex coefficients and
multiple eigenvalues are considered in \cite{BondarenkoGaidel2021}.
Transmutation operators for quadratic pencils are considered in \cite%
{AmirovGulyaz2014,Guseinov2014}.

To derive the NSBF representations for solutions of (\ref{eq:pencil}) we
expand the two associated transmutation operator kernels in
Fourier--Legendre series, thus applying the approach from \cite{KNT2017}.
Substitution of the Fourier--Legendre series into the integral
representations for the solutions gives their NSBF expansions. The same
approach applies to the integral representations of their first derivatives.
The resulting NSBF representations retain useful properties of those for the
Schr\"{o}dinger equation: they admit truncation bounds on horizontal strips
of the spectral plane that are independent of $|\operatorname{Re}\rho |$, and the
coefficients $p$ and $q$ can be recovered locally from the first NSBF
coefficients. These two features make the obtained NSBF representations
convenient for solving both direct and inverse problems for (\ref{eq:pencil}%
). Moreover, we develop a recursive integration procedure for computing the
NSBF coefficients when $p$ and $q$ are known. In a boundary value problem,
the fundamental solutions of (\ref{eq:pencil}) and their derivatives are
evaluated at many values and on large ranges of the spectral parameter $\rho
$. In the obtained NSBF representations, the coefficients depend only on $x$%
. The same computed coefficients can therefore be used throughout a
computation of a spectrum and for different separated boundary conditions.
The strip estimates control the truncation error as $|\operatorname{Re}\rho |$
grows. Together with the derivative representations, they also give error
bounds for characteristic functions of spectral problems and permit accurate
approximation of large sets of their isolated zeros.

The NSBF method for Schr\"{o}dinger and Sturm--Liouville equations was
developed in \cite{KNT2017,KT2017,KT2018,KravchenkoTorbaMMAS2021}.
Applications to inverse Sturm--Liouville problems with complex coefficients
are developed in \cite{Kravchenko2025Complex}. Spectral parameter power
series provide a related construction for polynomial pencils \cite%
{KravchenkoPorter2010,KravchenkoTorbaVelasco2015}. Further background can be
found in \cite%
{Carroll,Marchenko,Levitan1987,LevitanSargsjan1991,KravchenkoBook2020}.

The paper is organized as follows. Section~\ref{sec:transmutation} presents
the transmutation operators, the characteristic values of their kernels, and
the regularity needed below. In Section~\ref{sec:solution-nsbf} we derive
the NSBF representations for the solutions, and in Section~\ref%
{sec:truncation-estimates} we establish convergence and truncation
estimates. Section~\ref{sec:coefficient-construction} gives the recursive
integration procedure for the coefficients. Section~\ref%
{sec:potential-recovery} explains how the first coefficients can be used to
recover $p$ and $q$ directly. The representations for the first derivatives
and the corresponding coefficient formulas are obtained in Section~\ref%
{der:sec:derivative-nsbf}. Section~\ref{sec:computation} summarizes the
computational steps, while Section~\ref{sec:spectral-approximation} applies
the representations to separated boundary value problems. Two examples, an exactly
solvable one, and another with the known characteristic function, are considered in
Section~\ref{sec:example}, followed by the conclusions. The proof of the kernel
regularity lemma is given in the Appendix.

\section{Transmutation operators and integral representations}

\label{sec:transmutation}

Define
\begin{equation}
\sigma (x)=\int_{0}^{x}p(t)\,dt.  \label{eq:phase}
\end{equation}%
Denote
\begin{equation*}
\vartheta(\rho,x)=\rho x-\sigma(x).
\end{equation*}
Denote by $c$ and $s$ the solutions of the equation (\ref{eq:pencil})
defined by the initial conditions
\begin{equation}
\begin{split}
c(\rho ,0)& =1,\qquad c^{\prime }(\rho ,0)=0, \\
s(\rho ,0)& =0,\qquad s^{\prime }(\rho ,0)=1.
\end{split}
\label{rec:eq:c-s-initial}
\end{equation}
Additionally, consider the solutions
\begin{equation*}
e_{\pm }(\rho ,x)=c(\rho ,x)\pm i\rho s(\rho ,x).
\end{equation*}%
Their initial values take the form
\begin{equation}
e_{\pm }(\rho ,0)=1,\qquad e_{\pm }^{\prime }(\rho ,0)=\pm i\rho .
\label{rec:eq:e-initial-data}
\end{equation}

Define
\begin{equation*}
\mathcal{R}(x)=\int_{0}^{x}\left( q(t)+p(t)^{2}\right) \,dt.
\end{equation*}%
We begin by introducing integral representations for the solutions $e_{\pm
}(\rho ,x)$ in the form of Volterra integral operators applied to the
exponential functions $e^{\pm i\rho t}$. The Volterra operators of this kind
are known as transmutation operators and play an essential
role in spectral theory. Transmutation operators for real quadratic pencils
go back to \cite{GasymovGuseinov1981,Guseinov1986,Guseinov2014}.
Integral representations, characteristic data, and hyperbolic kernel
equations for the auxiliary solutions $e_\pm$ with real-valued
coefficients appear in \cite[Theorems~2.1 and~3.1, formulas~(35)--(36)
and~(53)]{AmirovGulyaz2014}. The following theorem gives a direct
construction for complex-valued $p$ and $q$. The approach is standard
and follows the reasoning in the preceding papers.

\begin{theorem}[Transmutation kernels and characteristic data]
\label{thm:transmutation-representation} Assume $q\in L^2(0,b)$ and $p\in
W^{1,2}(0,b)$. There exist unique continuous kernels $K_\pm(x,t)$ on $%
0\leq x\leq b$, $|t|\leq x$, such that, for every $x\in[0,b]$ and $\rho\in%
\mathbb{C}$,
\begin{equation}
e_\pm(\rho,x) =e^{\pm i\vartheta(\rho,x)} +\int _{-x}^{x}K_\pm(x,t)e^{\pm
i\rho t}\,dt,  \label{rec:eq:transmutation-representation}
\end{equation}
where $\sigma$ is defined by (\ref{eq:phase}). The characteristic data of
the kernels are
\begin{equation}
K_\pm(x,x) =\frac{e^{\mp i\sigma(x)}}2 \left(\mathcal{R}(x)\pm ip(x)\right),
\label{rec:eq:kernel-upper-value}
\end{equation}

\begin{equation}
K_{\pm }(x,-x)=\pm \frac{i}{2}p(0)e^{\pm i\sigma (x)}.
\label{rec:eq:kernel-lower-value}
\end{equation}%
If $q\in C^{1}[0,b]$ and $p\in C^{2}[0,b]$, each kernel is a classical
solution of the equation
\begin{equation}
K_{\pm ,xx}-K_{\pm ,tt}\mp 2ip(x)K_{\pm ,t}-q(x)K_{\pm }=0.
\label{rec:eq:kernel-pde}
\end{equation}
\end{theorem}

\noindent \textit{Proof.}

The hyperbolic equation (\ref{rec:eq:kernel-pde}) appears in
\cite[Theorem~3.1, formula~(53)]{AmirovGulyaz2014}; see also
\cite[Theorem~2.1]{Guseinov2014}. Denote by $\overline\Omega_b$ the closed triangle $\{(\xi,\eta):\xi\geq0,\
\eta\geq0,\ \xi+\eta\leq b\}$, and by $\Omega_b$ its interior. With
\begin{equation*}
\xi=(x+t)/2,\qquad \eta=(x-t)/2,\qquad
H_\pm(\xi,\eta)=K_\pm(\xi+\eta,\xi-\eta),
\end{equation*}
equation (\ref{rec:eq:kernel-pde}) takes the form
\begin{equation}
H_{\pm,\xi\eta}=\pm ip(\xi+\eta) (H_{\pm,\xi}-H_{\pm,\eta})+q(\xi+\eta)H_\pm.
\label{eq:goursat-H}
\end{equation}
Denote the right sides of (\ref{rec:eq:kernel-upper-value}) and (\ref%
{rec:eq:kernel-lower-value}) by
\begin{equation}
\begin{split}
h_{\pm,1}(\xi)&=\tfrac12e^{\mp i\sigma(\xi)} (\mathcal{R}(\xi)\pm ip(\xi)),
\\
h_{\pm,2}(\eta)&=\pm\tfrac i2p(0)e^{\pm i\sigma(\eta)}.
\end{split}
\label{eq:goursat-boundary-functions}
\end{equation}
They agree at the origin, and the corresponding characteristic conditions
take the form
\begin{equation}
H_\pm(\xi,0)=h_{\pm,1}(\xi),\qquad H_\pm(0,\eta)=h_{\pm,2}(\eta).
\label{eq:goursat-data}
\end{equation}
Integration of (\ref{eq:goursat-H}) over the rectangle $[0,\xi]\times[0,\eta]$,
followed by integration by parts, gives the equality
\begin{equation}
\begin{split}
H_\pm(\xi,\eta)={}&h_{\pm,1}(\xi)+h_{\pm,2}(\eta)-h_{\pm,1}(0) \\
&\pm i\int_0^\eta p(\xi+v)H_\pm(\xi,v)\,dv \mp i\int_0^\xi
p(u+\eta)H_\pm(u,\eta)\,du \\
&\pm i\int_0^\xi p(u)h_{\pm,1}(u)\,du \mp i\int_0^\eta p(v)h_{\pm,2}(v)\,dv
\\
&+\int_0^\xi\int_0^\eta q(u+v)H_\pm(u,v)\,dv\,du.
\end{split}
\label{eq:goursat-volterra}
\end{equation}
For $U\in C(\overline\Omega_b)$, define the integral operator
\begin{equation*}
\begin{split}
(\mathcal{V}_\pm U)(\xi,\eta)={}& \pm i\int_0^\eta p(\xi+v)U(\xi,v)\,dv \mp
i\int_0^\xi p(u+\eta)U(u,\eta)\,du \\
&+\int_0^\xi\int_0^\eta q(u+v)U(u,v)\,dv\,du.
\end{split}%
\end{equation*}
These are the three terms involving $H_\pm$ in (\ref{eq:goursat-volterra}),
with $H_\pm$ replaced by $U$. For $0\leq x\leq b$, let
\begin{equation*}
m_U(x)=\max_{\substack{ \xi,\eta\geq0  \\ \xi+\eta\leq x}}|U(\xi,\eta)|,
\qquad A(x)=\int_0^x(2|p(s)|+b|q(s)|)\,ds.
\end{equation*}
Thus $m_{\mathcal{V}_\pm U}(x)$ denotes the maximum of $|(\mathcal{V}_\pm
U)(\xi,\eta)|$ over the same triangle. It follows that
\begin{equation*}
m_{\mathcal{V}_\pm U}(x)\leq \int_0^x(2|p(s)|+b|q(s)|)m_U(s)\,ds.
\end{equation*}
Iteration therefore leads to the estimate
\begin{equation*}
m_{\mathcal{V}_\pm^nU}(x)\leq m_U(x)\frac{A(x)^n}{n!}.
\end{equation*}
Since $p,q\in L^1(0,b)$, the successive approximations converge uniformly on
$\overline\Omega_b$ and give the unique continuous solution of (\ref%
{eq:goursat-volterra}). If two sets of coefficients and characteristic data
are given, their solutions satisfy the difference of two equations of the
form (\ref{eq:goursat-volterra}). The same estimate bounds the difference of
these solutions by the differences of their data and by the $L^1$
differences of $p$ and $q$, provided the $L^1$ norms of both coefficient
pairs remain bounded.

For smooth coefficients the integral equation yields a $C^{2}$ kernel.
Differentiation of the characteristic data (\ref%
{eq:goursat-boundary-functions}) shows that its derivatives satisfy the
equalities
\begin{equation*}
\begin{split}
(K_{\pm ,x}+K_{\pm ,t}\pm ipK_{\pm })(x,x)&=\tfrac{1}{2}(q(x)+p^{2}(x)\pm
ip^{\prime }(x))e^{\mp i\sigma (x)}, \\
(K_{\pm ,x}-K_{\pm ,t}\mp ipK_{\pm })(x,-x)&=0.
\end{split}%
\end{equation*}%
Integrating twice by parts in (\ref{rec:eq:transmutation-representation}),
using (\ref{rec:eq:kernel-pde}) and these identities, proves (\ref{eq:pencil}%
) and the initial values (\ref{rec:eq:e-initial-data}). Uniqueness of the
initial value problem proves (\ref{rec:eq:transmutation-representation}) in
the smooth case. Finally, take smooth $q_{\nu }\rightarrow q$ in $L^{2}$ and
$p_{\nu }\rightarrow p$ in $W^{1,2}$. The characteristic data converge
uniformly because $p_{\nu }\rightarrow p$ in $C[0,b]$, while $q_{\nu
}\rightarrow q$ in $L^{1}$. The preceding Volterra estimate gives uniform
convergence of the kernels. Continuous dependence of the initial value
solutions, written as first order integral systems, gives uniform
convergence of the solutions and their first $x$ derivatives. Passing to the
limit establishes the representation. Uniqueness of a continuous kernel in
that representation also follows from injectivity of the Fourier transform
of its restriction to $[-x,x]$ for each $x$. \hfill $\square $

The kernels of Theorem~\ref{thm:transmutation-representation} define
operators $\mathbf{T}_\pm:C[-b,b]\to C[0,b]$ by
\begin{equation}
\mathbf{T}_\pm[v](x) =e^{\mp i\sigma(x)}v(x)
+\int_{-x}^{x}K_\pm(x,t)v(t)\,dt.  \label{rec:eq:transmutation-operator}
\end{equation}
By (\ref{rec:eq:transmutation-operator}),
\begin{equation*}
e_\pm(\rho,x) =\mathbf{T}_\pm\left[e^{\pm i\rho\,\cdot}\right](x).
\end{equation*}

Theorem~2.1 of \cite{Guseinov2014} states two-variable Sobolev regularity
for the corresponding sine and cosine kernels with real coefficients. In the next lemma,
we obtain $L^2$ estimates for $K_\pm(x,\cdot)$ and their derivatives for
complex-valued coefficients $p$ and $q$, uniformly for $0<x\leq b$.

\begin{lemma}[Regularity of the transmutation kernels]
\label{lem:goursat-sobolev-regularity} Let $m=0,1,\ldots$, $q\in
W^{m,2}(0,b)$, and $p\in W^{m+1,2}(0,b)$. The solution $H_\pm$ of (\ref%
{eq:goursat-volterra}) belongs to $W^{m+1,2}(\Omega_b)$ and satisfies (\ref%
{eq:goursat-H}) in the sense of distributions in $\Omega_b$, with the traces (\ref%
{eq:goursat-data}). Moreover,
\begin{align}
\sup_{0<x\leq b}\sum_{j=0}^{m+1} \|\partial_t^j
K_\pm(x,\cdot)\|_{L^2(-x,x)}^2 &\leq C_{m,b,q,p}^{\,2},\label{eq:goursat-section-bound} \\
\sup_{0<x\leq b}\sum_{j=0}^{m} \|\partial_t^j
K_{\pm,x}(x,\cdot)\|_{L^2(-x,x)}^2 &\leq C_{m,b,q,p}^{\,2}.\label{eq:goursat-section-bound-derivative}
\end{align}
Here $C_{m,b,q,p}$ is locally bounded as a function of $\|q\|_{W^{m,2}}$ and
$\|p\|_{W^{m+1,2}}$, and the same constant may be used for both signs. The $%
t $ derivatives in (\ref{eq:goursat-section-bound}), (\ref{eq:goursat-section-bound-derivative}) exist for every $%
x\in(0,b]$. After extension by zero from $(-x,x)$ to $(-b,b)$, each function
$x\mapsto\partial_t^jK_\pm(x,\cdot)$ and $x\mapsto\partial_t^jK_{\pm,x}(x,%
\cdot)$ appearing there is continuous with values in $L^2(-b,b)$, with value
zero at $x=0$.
\end{lemma}

Analogous regularity results are stated in \cite[Theorem~2.1]{Guseinov2014} for real-valued $p$ and $q$,
and in \cite[Propositions~4.1 and~4.7]{RoqueTorba2025} for Dirac systems. We give the details for both complex kernels in the Appendix for completeness.

\section{Neumann series of Bessel functions for solutions}

\label{sec:solution-nsbf}

Let $P_{n}$ denote the Legendre polynomial of degree $n$. For $x>0$,
introduce the Fourier--Legendre coefficients
\begin{equation*}
\kappa _{n}^{\pm }(x)=\frac{2n+1}{2}\int_{-x}^{x}K_{\pm }(x,t)P_{n}\left(
\frac{t}{x}\right) \,dt.
\end{equation*}%
Then
\begin{equation}
K_{\pm }(x,t)=\sum_{n=0}^{\infty }\frac{\kappa _{n}^{\pm }(x)}{x}P_{n}\left(
\frac{t}{x}\right) .  \label{rec:eq:kernel-legendre-expansion}
\end{equation}%
The equality holds in $L^{2}(-x,x)$.

The continuity of the kernels gives
\[
|\kappa_n^\pm(x)|\leq (2n+1)x\max_{|t|\leq x}|K_\pm(x,t)|\longrightarrow0,\quad x\to 0,
\]
hence the coefficients $\kappa _{n}^{\pm }$ can be continuously extended to 0 with the limit values
\[
\kappa_n^\pm(0)=0,\qquad n=0,1,2,\ldots.
\]

We now convert the kernel expansions into NSBF representations. The required
Legendre and Bessel identities can be found in \cite{AbramowitzStegun}, and
the approximation statements rely on the standard theory summarized in \cite%
{DeVoreLorentz}. We write $J_{\nu }$ for the Bessel function of the first
kind of order $\nu $. For $n\in \mathbb{N}_{0}$ and $z\in \mathbb{C}$, the
spherical Bessel function of order $n$ is denoted by $j_{n}$ and defined by
\begin{equation*}
j_{n}(z)=\sqrt{\frac{\pi }{2z}}J_{n+1/2}(z).
\end{equation*}%
In classical analysis, a Neumann series of Bessel functions is a series of
the form
\begin{equation*}
\sum_{n=0}^{\infty}a_nJ_{\nu+n}(z),
\end{equation*}
where $\nu$ is fixed and the coefficients $a_n$ are independent of $z$. The
Bessel functions have a common argument and their orders increase with the
summation index. The basic theory can be found in \cite[p.~522]{Watson}
and \cite{Wilkins1948}. The relation between $j_n$ and $J_{n+1/2}$ explains
the use of the same terminology for the spherical Bessel series and their
even and odd subseries considered here.

For each fixed $n$, the small argument expansions give
\begin{equation*}
j_0(z)=1+O(z^2),\qquad j_1(z)=\frac{z}{3}+O(z^3),\qquad
j_n(z)=O(z^n)\quad(n\geq2).
\end{equation*}%
In particular, $j_0(0)=1$ and $j_n(0)=0$ for $n\geq1$ \cite[Eq.~(10.53.1)]%
{DLMF}. We use the integral identity (the Fourier transform of the Legendre
polynomials)
\begin{equation}
\int_{-1}^{1}P_{n}(t)e^{izt}\,dt=2i^{n}j_{n}(z).
\label{eq:legendre-bessel-integral}
\end{equation}

\begin{theorem}[NSBF representations]
\label{Th NSBF}For every $0<x\leq b$, $\rho \in \mathbb{C}$, the solutions $%
e_{\pm }(\rho ,x)$ admit the following series representations
\begin{equation}
e_{\pm }(\rho ,x)=e^{\pm i\vartheta (\rho ,x)}+2\sum_{n=0}^{\infty }(\pm
i)^{n}\kappa _{n}^{\pm }(x)j_{n}(\rho x).  \label{rec:eq:e-pm-nsbf}
\end{equation}%
For each fixed $x>0$ and $R>0$, the series converge absolutely and uniformly
for $|\rho |\leq R$.
\end{theorem}

\noindent \textit{Proof.} Denote the Fourier--Legendre partial sum by
\begin{equation*}
K_{\pm ,N}(x,t)=\sum_{n=0}^{N}\frac{\kappa _{n}^{\pm }(x)}{x}P_{n}(t/x).
\end{equation*}%
First replace $K_{\pm }$ in the transmutation integral by $K_{\pm ,N}$. With
$t=xu$, the identity (\ref{eq:legendre-bessel-integral}) gives
\begin{equation*}
\int_{-x}^{x}K_{\pm ,N}(x,t)e^{\pm i\rho t}\,dt =\sum_{n=0}^{N}\kappa _{n}^{\pm }(x)\int_{-1}^{1}P_{n}(u)e^{\pm
i\rho xu}\,du=2\sum_{n=0}^{N}(\pm i)^{n}\kappa _{n}^{\pm }(x)j_{n}(\rho x).
\end{equation*}%
The partial sums converge to $K_{\pm }(x,\cdot )$ in $L^{2}(-x,x)$. For $%
|\rho |\leq R$, the $L^{2}(-x,x)$ norm of $e^{\pm i\rho t}$ does not exceed $%
\sqrt{2x}e^{Rx}$. The Cauchy--Schwarz inequality therefore permits passage
to the limit uniformly on this disk and proves (\ref{rec:eq:e-pm-nsbf}).

To verify absolute convergence, expand the exponential in (\ref%
{eq:legendre-bessel-integral}) into its Taylor series. The polynomial $P_{n}$
is orthogonal to every polynomial of degree less than $n$. The Taylor
remainder and $\int_{-1}^{1}|P_{n}(t)|\,dt\leq 2$ imply
\begin{equation*}
|j_{n}(\rho x)|\leq e^{Rx}\frac{(Rx)^{n}}{n!},\qquad |\rho |\leq R.
\end{equation*}%
The definition of $\kappa _{n}^{\pm }$ and the Cauchy--Schwarz inequality
give
\begin{equation*}
|\kappa _{n}^{\pm }(x)|\leq \left( \frac{x(2n+1)}{2}\right) ^{1/2}\Vert
K_{\pm }(x,\cdot )\Vert _{L^{2}(-x,x)}.
\end{equation*}%
Since $\sum_{n=0}^{\infty }\sqrt{2n+1}(Rx)^{n}/n!$ converges, these
estimates prove absolute and uniform convergence on $|\rho |\leq R$ (cf.
\cite{KNT2017}). \hfill $\square $

Define the four coefficient families associated with the solutions $c(\rho
,x)$ and $s(\rho ,x)$\ by
\begin{align}
a_{n}(x)&=\kappa _{2n}^{+}(x)+\kappa _{2n}^{-}(x),& b_{n}(x)&=i\left(
\kappa _{2n+1}^{+}(x)-\kappa _{2n+1}^{-}(x)\right) ,
\label{rec:eq:a-b-definitions}\\
u_{n}(x)&=-i\left( \kappa _{2n}^{+}(x)-\kappa _{2n}^{-}(x)\right) ,&
v_{n}(x)&=\kappa _{2n+1}^{+}(x)+\kappa _{2n+1}^{-}(x).
\label{rec:eq:u-v-definitions}
\end{align}%
Then (\ref{rec:eq:a-b-definitions}) and (\ref{rec:eq:u-v-definitions}) give
the following corollary of Theorem~\ref{Th NSBF}.

\begin{corollary}
The solutions of (\ref{eq:pencil}) specified by (\ref{rec:eq:c-s-initial})
admit the series representations
\begin{align}
c(\rho ,x)&= \cos \vartheta (\rho ,x) +\sum_{n=0}^{\infty}(-1)^{n}a_{n}(x)j_{2n}(\rho x) +\sum_{n=0}^{\infty }(-1)^{n}b_{n}(x)j_{2n+1}(\rho x),
\label{rec:eq:c-main}\\
s(\rho ,x)&= \frac{\sin \vartheta (\rho ,x)+u_{0}(x)j_{0}(\rho x)}{\rho }
 +\frac{1}{\rho }\sum_{n=1}^{\infty }(-1)^{n}u_{n}(x)j_{2n}(\rho x)
 +\frac{1}{\rho }\sum_{n=0}^{\infty }(-1)^{n}v_{n}(x)j_{2n+1}(\rho x).
\label{rec:eq:s-main}
\end{align}%
For each fixed $x>0$ and $R>0$, the series in (\ref{rec:eq:c-main}) and (\ref%
{rec:eq:s-main}) converge absolutely and uniformly for $|\rho |\leq R$. The small argument expansions of the spherical Bessel functions lead to the expressions
\begin{equation*}
c(0,x)=\cos \sigma (x)+a_{0}(x),\qquad s(0,x)=x\cos \sigma (x)+\frac{x}{3}%
v_{0}(x).
\end{equation*}%
At $x=0$, the formulas give $c(\rho,0)=1$ and $s(\rho,0)=0$.
\end{corollary}

\section{Convergence and truncation estimates}

\label{sec:truncation-estimates}

For computations concerning direct and inverse spectral problems, one also
needs to know how accurately a finite number of terms represents the
solutions when $|\operatorname{Re}\rho |$ is large. We now estimate the truncation
errors from the Fourier--Legendre remainders of the kernels. The estimates
hold on entire horizontal strips and give rates in the truncation index when
the coefficients in the pencil have additional regularity.

For $x>0$ and $N\geq 0$, let
\begin{equation*}
R_{\pm ,N}(x,t)=K_{\pm }(x,t)-K_{\pm ,N}(x,t),
\end{equation*}%
and define
\begin{equation*}
d_{\pm ,N}(x)=\Vert R_{\pm ,N}(x,\cdot )\Vert _{L^{2}(-x,x)}.
\end{equation*}%
The orthogonality relations for the Legendre polynomials and Parseval's
identity give
\begin{equation}
d_{\pm ,N}(x)^{2}=\frac{2}{x}\sum_{n=N+1}^{\infty }\frac{|\kappa _{n}^{\pm
}(x)|^{2}}{2n+1}.  \label{err:eq:parseval-tail}
\end{equation}%
Let $e_{\pm ,N}(\rho ,x)$ be obtained by retaining the terms with $0\leq
n\leq N$ in (\ref{rec:eq:e-pm-nsbf}), and define
\begin{equation*}
c_{N}(\rho ,x)=\frac{e_{+,N}(\rho ,x)+e_{-,N}(\rho ,x)}{2},\qquad s_{N}(\rho
,x)=\frac{e_{+,N}(\rho ,x)-e_{-,N}(\rho ,x)}{2i\rho },
\end{equation*}%
At $\rho =0$, these truncations have the exact values
\begin{equation*}
e_{\pm ,N}(0,x)=e_{\pm }(0,x)=c_{N}(0,x)=c(0,x),\qquad s_{N}(0,x)=%
\begin{cases}
x\cos \sigma (x), & N=0, \\
x\cos \sigma (x)+xv_{0}(x)/3, & N\geq 1.%
\end{cases}%
\end{equation*}%
After division by $\rho $, only the term involving $j_{1}$ contributes at
zero, and it is retained when $N\geq 1$.

We use the two auxiliary functions
\begin{equation}
\Lambda _{C}(x)=%
\begin{cases}
\left( \frac{\sinh (2Cx)}{C}\right) ^{1/2}, & C>0, \\[2mm]
\sqrt{2x}, & C=0,%
\end{cases}
\label{err:eq:Lambda}
\end{equation}%
and
\begin{equation}
\Theta _{C}(x)=%
\begin{cases}
\left( \frac{e^{2Cx}(2C^{2}x^{2}-2Cx+1)-1}{2C^{3}}\right) ^{1/2}, & C>0, \\%
[2mm]
\sqrt{\frac{2x^{3}}{3}}, & C=0.%
\end{cases}
\label{err:eq:Theta}
\end{equation}%
Thus, for $C\geq 0$,
\begin{equation*}
\sup_{|\operatorname{Im}\rho |\leq C}\int_{-x}^{x}|e^{\pm i\rho t}|^{2}\,dt=\Lambda
_{C}(x)^{2},\qquad \Theta _{C}(x)^{2}=\int_{-x}^{x}t^{2}e^{2C|t|}\,dt,
\end{equation*}%
and $\Lambda _{C}(x)\leq \sqrt{2x}e^{Cx}$, $\Theta _{C}(x)\leq e^{Cx}\sqrt{%
2x^{3}/3}$.

For $z\in \mathbb{C}$, the following identity is valid
\begin{equation}
\sum_{n=0}^{\infty }(2n+1)|j_{n}(z)|^{2}=%
\begin{cases}
\frac{\sinh (2|\operatorname{Im}z|)}{2|\operatorname{Im}z|}, & \operatorname{Im}z\neq 0, \\[2mm]
1, & \operatorname{Im}z=0.%
\end{cases}%
.  \label{err:eq:bessel-parseval}
\end{equation}

It follows from the spherical Bessel addition formula \cite[Eq.~(10.60.2)]%
{DLMF}.

\begin{theorem}[Uniform NSBF error estimates on horizontal strips]
\label{thm:uniform-error-estimates} Let $N\geq 0$, $0<x\leq b$, $C\geq 0$,
and $|\operatorname{Im}\rho |\leq C$. Then
\begin{equation}
|e_{\pm }(\rho ,x)-e_{\pm ,N}(\rho ,x)|\leq \Lambda _{C}(x)d_{\pm ,N}(x),
\label{err:eq:e-estimate}
\end{equation}%
\begin{equation}
|c(\rho ,x)-c_{N}(\rho ,x)|\leq \frac{\Lambda _{C}(x)}{2}\left(
d_{+,N}(x)+d_{-,N}(x)\right) ,  \label{err:eq:c-estimate}
\end{equation}%
and
\begin{equation}
|\rho |\,|s(\rho ,x)-s_{N}(\rho ,x)|\leq \frac{\Lambda _{C}(x)}{2}\left(
d_{+,N}(x)+d_{-,N}(x)\right) .  \label{err:eq:rho-s-estimate}
\end{equation}%
For the error in $s$ without the factor $|\rho |$, we have
\begin{equation}
|s(\rho ,x)-s_{N}(\rho ,x)|\leq \frac{\Theta _{C}(x)}{2}\left(
d_{+,N}(x)+d_{-,N}(x)\right) .  \label{err:eq:s-estimate}
\end{equation}%
At $\rho =0$, these errors have the exact values
\begin{equation*}
e_{\pm }(0,x)-e_{\pm ,N}(0,x)=c(0,x)-c_{N}(0,x)=0,\qquad s(0,x)-s_{N}(0,x)=%
\begin{cases}
xv_{0}(x)/3, & N=0, \\
0, & N\geq 1.%
\end{cases}%
\end{equation*}%
For each fixed $x>0$, the series in (\ref{rec:eq:e-pm-nsbf}) and (\ref%
{rec:eq:c-main}), as well as the two series in the regularized formula (\ref%
{rec:eq:s-main}), converge absolutely and uniformly on the whole strip $|%
\operatorname{Im}\rho |\leq C$. None of these bounds grows with $|\operatorname{Re}\rho |$.
\end{theorem}

\noindent \textit{Proof.} By (\ref{err:eq:parseval-tail}), (\ref%
{err:eq:bessel-parseval}), (\ref{err:eq:Lambda}), and the Cauchy--Schwarz
inequality, we have
\begin{equation*}
2\sum_{n>N}|\kappa _{n}^{\pm }(x)j_{n}(\rho x)|\leq \Lambda _{C}(x)d_{\pm
,N}(x).
\end{equation*}%
Since $d_{\pm ,N}(x)\rightarrow 0$, the exponential solution series
converges absolutely and uniformly on the whole strip. The same inequality
proves (\ref{err:eq:e-estimate}). This yields (\ref{err:eq:c-estimate}) and (%
\ref{err:eq:rho-s-estimate}) as well.

Consider $s-s_{N}$ without multiplying by $|\rho |$. For $n>N\geq 0$, $P_{n}$
is orthogonal to constants. For $\rho \neq 0$, its Legendre coefficient in $%
e^{\pm i\rho t}/\rho $ agrees with that in $(e^{\pm i\rho t}-1)/\rho $. The
latter function extends through zero because
\begin{equation*}
\frac{e^{\pm i\rho t}-1}{\rho }=\pm it\int_{0}^{1}e^{\pm i\theta \rho
t}\,d\theta .
\end{equation*}%
On $|\operatorname{Im}\rho |\leq C$ its modulus does not exceed $|t|e^{C|t|}$, and
its $L^{2}(-x,x)$ norm does not exceed $\Theta _{C}(x)$ by (\ref%
{err:eq:Theta}). Parseval's identity and the Cauchy--Schwarz inequality, applied to these
Legendre coefficients, yield
\begin{equation*}
2\sum_{n>N}\left\vert \kappa _{n}^{\pm }(x)\frac{j_{n}(\rho x)}{\rho }%
\right\vert \leq \Theta _{C}(x)d_{\pm ,N}(x),
\end{equation*}%
with their values at $\rho =0$ given by the small argument expansions of $%
j_{n}$. Taking half the difference proves (\ref{err:eq:s-estimate}). For $%
N=0 $, the omitted $j_{1}$ term contributes $xv_{0}(x)/3$ at zero. For $%
N\geq 1$, every omitted term vanishes after division by $\rho $. Since $%
d_{\pm ,N}(x)\rightarrow 0$, this also proves absolute uniform convergence
of the two divided series in (\ref{rec:eq:s-main}) on the strip. The first
absolute tail bound gives the same conclusion for both series in (\ref%
{rec:eq:c-main}). \hfill $\square $

We combine the Goursat regularity result with an estimate for Fourier--Legendre remainders to obtain bounds in terms of $q$ and $p$.

\begin{theorem}[Coefficient-based Fourier--Legendre rates]
\label{thm:kernel-projection-rates} Let $r\in\mathbb{N}$ and assume
\begin{equation*}
q\in W^{r-1,2}(0,b),\qquad p\in W^{r,2}(0,b).
\end{equation*}
For every $0<x\leq b$, the kernels $K_\pm(x,\cdot)$ belong to $W^{r,2}(-x,x)$%
. Define
\begin{equation}
\begin{split}
\mathcal{A}_r^\pm&:=\sup_{0<x\leq b} \|\partial_t^r
K_\pm(x,\cdot)\|_{L^2(-x,x)}, \\
\mathcal{A}_r&:=\max\{\mathcal{A}_r^+,\mathcal{A}_r^-\}<\infty.
\end{split}
\label{err:eq:coefficient-kernel-bound}
\end{equation}
For fixed $b$, Lemma~\ref{lem:goursat-sobolev-regularity} bounds $\mathcal{A}%
_r$ by a locally bounded function of $\|q\|_{W^{r-1,2}}$ and $%
\|p\|_{W^{r,2}} $. In all the estimates below, the same positive constant $%
C_r$ can be used for both signs. It depends only on $r$ and is independent
of the coefficients, $b$, $x$, $N$, $\rho$ and the strip width $C$.

For $0<x\leq b$ and $N>r$, we have
\begin{equation}
d_{\pm,N}(x)\leq\frac{C_r x^r}{N^r} \left\|\left(1-\frac{t^2}{x^2}%
\right)^{r/2} \partial_t^r K_\pm(x,t)\right\|_{L^2(-x,x)}.
\label{err:eq:kernel-L2-rate}
\end{equation}
In particular,
\begin{equation}
d_{\pm,N}(x)\leq\frac{C_r\mathcal{A}_r x^r}{N^r}.
\label{err:eq:coefficient-L2-rate}
\end{equation}
For $n>r+1$, the kernel coefficients satisfy
\begin{equation}
|\kappa_n^\pm(x)|\leq C_r\mathcal{A}_r x^{r+1/2}n^{-r+1/2}.
\label{err:eq:kappa-decay}
\end{equation}

For $0<x\leq b$, $C\geq 0$, $|\operatorname{Im}\rho |\leq C$ and $N>r$, we have
\begin{equation}
|e_{\pm }(\rho ,x)-e_{\pm ,N}(\rho ,x)|\leq \frac{C_{r}\mathcal{A}%
_{r}x^{r}\Lambda _{C}(x)}{N^{r}},  \label{err:eq:e-sobolev-rate}
\end{equation}%
\begin{equation}
\max \left\{ |c(\rho ,x)-c_{N}(\rho ,x)|,\,|\rho |\,|s(\rho ,x)-s_{N}(\rho
,x)|\right\} \leq \frac{C_{r}\mathcal{A}_{r}x^{r}\Lambda _{C}(x)}{N^{r}}.
\label{err:eq:fundamental-sobolev-rate}
\end{equation}%
The estimate for $s$ without the factor $|\rho |$ can be written as
\begin{equation}
|s(\rho ,x)-s_{N}(\rho ,x)|\leq \frac{C_{r}\mathcal{A}_{r}x^{r}\Theta _{C}(x)%
}{N^{r}}.  \label{err:eq:s-sobolev-rate}
\end{equation}
\end{theorem}

\noindent\textit{Proof.} Apply Lemma~\ref{lem:goursat-sobolev-regularity}
with $m=r-1$. It gives the stated Sobolev regularity and the uniform bound (%
\ref{err:eq:coefficient-kernel-bound}).

Fix $x>0$ and one sign, and set $F(u)=K_\pm(x,xu)$ and $F_N(u)=K_{%
\pm,N}(x,xu)$ for $-1<u<1$. Then $F_N$ is the $N$th partial sum of the Fourier--Legendre series of $F$.

The orthogonality relations for associated Legendre polynomials \cite[Sec.~8.14]{AbramowitzStegun} give,
first for polynomials and then by density for complex valued $F\in
W^{r,2}(-1,1)$, the following estimate for the remainder (see also \cite{Ky1992}):
\begin{equation*}
\Vert F-F_{N}\Vert _{L^{2}(-1,1)}\leq \sqrt{\frac{(N+1-r)!}{(N+1+r)!}}\Vert
(1-u^{2})^{r/2}F^{(r)}\Vert _{L^{2}(-1,1)}.
\end{equation*}%
The change of variables $t=xu$ gives
\begin{equation*}
F^{(r)}(u)=x^{r}\partial _{t}^{r}K_{\pm }(x,xu),\qquad d_{\pm ,N}(x)=\sqrt{x}%
\,\Vert F-F_{N}\Vert _{L^{2}(-1,1)}.
\end{equation*}%
Consequently, the bound before replacing the factorials by a power
of $N$ has the form
\begin{equation*}
d_{\pm ,N}(x)\leq x^{r}\sqrt{\frac{(N+1-r)!}{(N+1+r)!}}\left\Vert \left( 1-%
\frac{t^{2}}{x^{2}}\right) ^{r/2}\partial _{t}^{r}K_{\pm }(x,t)\right\Vert
_{L^{2}(-x,x)}.
\end{equation*}%
For fixed $r$, the factorial factor is bounded by $C_{r}N^{-r}$ for $N>r$.
This proves (\ref{err:eq:kernel-L2-rate}). The weight does not exceed one,
so (\ref{err:eq:coefficient-kernel-bound}) also gives (\ref%
{err:eq:coefficient-L2-rate}).

Because $P_n(t/x)$ is orthogonal to the polynomial $K_{\pm,n-1}(x,t)$, the
definition of $\kappa_n^\pm$ and the Cauchy--Schwarz inequality give
\begin{equation*}
|\kappa_n^\pm(x)|\leq \sqrt{\frac{x(2n+1)}2}\,d_{\pm,n-1}(x).
\end{equation*}
For $n>r+1$, apply the preceding factorial bound with $N=n-1$. Then
\begin{equation*}
|\kappa_n^\pm(x)|\leq \mathcal{A}_r x^{r+1/2} \sqrt{\frac{2n+1}{2}}\sqrt{%
\frac{(n-r)!}{(n+r)!}}.
\end{equation*}
The product of square roots is bounded by $C_r n^{-r+1/2}$. Choose $C_r$
large enough for both factorial bounds. This proves (\ref{err:eq:kappa-decay}%
) with the same constant. Finally, substitute (\ref%
{err:eq:coefficient-L2-rate}) into the strip estimates (\ref%
{err:eq:e-estimate})--(\ref{err:eq:s-estimate}). The resulting bounds are (%
\ref{err:eq:e-sobolev-rate}), (\ref{err:eq:fundamental-sobolev-rate}) and (%
\ref{err:eq:s-sobolev-rate}). \hfill$\square$

In particular, the hypotheses of Theorem~\ref%
{thm:transmutation-representation} give the case $r=1$. Since $x\leq b$, the resulting $N^{-1}$ estimates hold uniformly for $x\in[0,b]$.

Finally, bounded spectral sets admit a much faster estimate than the uniform
strip bound. No additional smoothness is needed.

\begin{proposition}[Factorial convergence on bounded spectral sets]
\label{prop:compact-factorial-convergence} Let $R>0$, $0<x\leq b$, $%
|\rho|\leq R$, and $N\geq0$. Then
\begin{equation}
|e_\pm-e_{\pm,N}| \leq d_{\pm,N}(x) \sqrt{\frac{2x}{2N+3}}\,e^{Rx} \frac{%
(Rx)^{N+1}}{(N+1)!}.  \label{err:eq:e-factorial}
\end{equation}
For the quotient of the error by $\rho$,
\begin{equation}
\left|\frac{e_\pm-e_{\pm,N}}{\rho}\right| \leq d_{\pm,N}(x) \sqrt{\frac{2x}{%
2N+3}}\,e^{Rx} \frac{x(Rx)^N}{(N+1)!}.  \label{err:eq:e-over-rho-factorial}
\end{equation}
The quotient on the left of (\ref{err:eq:e-over-rho-factorial}) has the value
\begin{equation*}
\left.\frac{e_\pm(\rho,x)-e_{\pm,N}(\rho,x)}{\rho}\right|_{\rho=0} =%
\begin{cases}
\pm 2ix\kappa_1^\pm(x)/3, & N=0, \\
0, & N\geq1.%
\end{cases}%
\end{equation*}%
The bound in (\ref{err:eq:e-over-rho-factorial}) holds for these values as
well. The half-sum and the quotient of the difference by $2i\rho$ give the
corresponding estimates for $c-c_N$ and $s-s_N$. Since $d_{\pm,N}\leq\|K_%
\pm(x,\cdot)\|_{L^2(-x,x)}$, these are factorial bounds in $N$ without an
additional smoothness assumption.
\end{proposition}

\noindent \textit{Proof.} The tail $R_{\pm ,N}$ is orthogonal to every
polynomial of degree at most $N$. Therefore, with $T_{N}(z)=%
\sum_{k=0}^{N}z^{k}/k!$,
\begin{equation*}
e_{\pm }-e_{\pm ,N}=\int_{-x}^{x}R_{\pm ,N}(x,t)\left( e^{\pm i\rho
t}-T_{N}(\pm i\rho t)\right) \,dt.
\end{equation*}%
Taylor's remainder estimate, followed by the Cauchy--Schwarz inequality,
gives the estimate
\begin{equation*}
|e_{\pm }-e_{\pm ,N}|\leq \frac{e^{Rx}|\rho |^{N+1}}{(N+1)!}d_{\pm
,N}(x)\left( \int_{-x}^{x}|t|^{2N+2}\,dt\right) ^{1/2},
\end{equation*}%
and
\begin{equation*}
\left( \int_{-x}^{x}|t|^{2N+2}\,dt\right) ^{1/2}=x^{N+1}\sqrt{\frac{2x}{2N+3}%
}.
\end{equation*}%
These inequalities prove (\ref{err:eq:e-factorial}). Division by $|\rho |$
before using $|\rho |\leq R$ proves (\ref{err:eq:e-over-rho-factorial}). The
tail is orthogonal to $1$, so its quotient by $\rho $ is regular at zero.
Since $j_{1}(z)=z/3+O(z^{3})$ and $j_{n}(z)=O(z^{n})$ for $n\geq 2$, the
first omitted coefficient gives the exact value stated above. \hfill $%
\square $

\section{Recursive construction of the coefficients}

\label{sec:coefficient-construction}

\subsection{Transmuted monomials and a direct coefficient formula}

\label{rec:sec:transmuted}

Denote the differential expression by
\begin{equation*}
L[y]=y^{\prime \prime }-q(x)y.
\end{equation*}

We first obtain a direct coefficient formula and then derive a recurrence
suitable for computation.

For $k\in\mathbb{N}_0$, define the corresponding transmuted monomial by
\begin{equation}
\Phi_{k}^\pm(x)=\mathbf{T}_\pm[t^{k}](x).  \label{rec:eq:Phi-definition}
\end{equation}
Expanding the exponential in (\ref{rec:eq:transmutation-representation}) and
integrating its power series term by term gives
\begin{equation}
e_\pm(\rho,x) =\sum_{k=0}^{\infty}\frac{(\pm i\rho)^{k}}{k!} \Phi_{k}^\pm(x).
\label{rec:eq:e-Taylor-Phi}
\end{equation}

\begin{proposition}
The functions $\Phi _{k}^{\pm }$ satisfy
\begin{equation}
L\left[ \Phi _{k}^{\pm }\right] =k(k-1)\Phi _{k-2}^{\pm }\mp 2ikp(x)\Phi
_{k-1}^{\pm },\qquad k=0,1,2,\ldots ,  \label{rec:eq:Phi-recurrence}
\end{equation}%
where $\Phi _{-1}^{\pm }=\Phi _{-2}^{\pm }=0$. Their initial values are
\begin{equation*}
\Phi _{k}^{\pm }(0)=\delta _{k0},\qquad (\Phi _{k}^{\pm })^{\prime
}(0)=\delta _{k1}\qquad (k\geq 0).
\end{equation*}%
Here $\delta _{jk}$ is the Kronecker symbol.
\end{proposition}

\noindent \textit{Proof.} The Volterra integral equation for the initial
value problem shows that $\rho \mapsto e_{\pm }(\rho ,\cdot )$ is entire
with values in $W^{2,2}(0,b)$. Its Taylor series can therefore be
differentiated twice with respect to $x$ in $L^{2}$. By (\ref%
{rec:eq:e-Taylor-Phi}), its coefficients are $\Phi _{k}^{\pm }$. Apply $L$
to the series and substitute it into
\begin{equation*}
L[e_{\pm }]+\rho ^{2}e_{\pm }-2\rho p(x)e_{\pm }=0.
\end{equation*}%
Set $z_{\pm }=\pm i\rho $. Then $\rho =\mp iz_{\pm }$ and $\rho ^{2}=-z_{\pm
}^{2}$. Comparison of the coefficient of $z_{\pm }^{k}/k!$ gives (\ref%
{rec:eq:Phi-recurrence}). The initial values follow by expanding $e_{\pm
}(\rho ,0)=1$ and $e_{\pm }^{\prime }(\rho ,0)=\pm i\rho $ in powers of $\pm
i\rho $. \hfill $\square $

Note that the functions $\Phi _{k}^{\pm }$ are closely related to the formal
powers introduced in \cite{KravchenkoTorbaVelasco2015}, where a recursive
integration procedure for their construction is developed.

Let
\begin{equation}
P_{n}(z)=\sum_{k=0}^{n}l_{k,n}z^{k}.  \label{rec:eq:Legendre-coefficients}
\end{equation}

\begin{proposition}
For every $n\geq 0$ and $x>0$ the equality is valid
\begin{equation}
\kappa _{n}^{\pm }(x)=\frac{2n+1}{2}\left( \sum_{k=0}^{n}l_{k,n}\frac{\Phi
_{k}^{\pm }(x)}{x^{k}}-e^{\mp i\sigma (x)}\right) .
\label{rec:eq:direct-kappa-formula}
\end{equation}%
\end{proposition}

\noindent\textit{Proof.} Using (\ref{rec:eq:Legendre-coefficients}), the
definition of the transmuted monomials gives
\begin{equation*}
\sum_{k=0}^n l_{k,n}\frac{\Phi_k^\pm(x)}{x^k} =\mathbf{T}_\pm[P_n(\cdot/x)](x)=e^{\mp i\sigma(x)}P_n(1) +\int_{-x}^xK_\pm(x,t)P_n(t/x)\,dt.
\end{equation*}
Since $P_n(1)=1$, the definition of $\kappa_n^\pm$ reduces this equality to (%
\ref{rec:eq:direct-kappa-formula}). \hfill$\square$

\begin{remark}
Formula (\ref{rec:eq:direct-kappa-formula}) is useful for symbolic checks,
but the recurrence below avoids the severe numerical cancellation caused by
the Legendre coefficients in this formula. See \cite{KNT2017}.
\end{remark}

\subsection{Differential recurrence for NSBF coefficients}

To derive the differential recurrence, assume temporarily that $q,p\in C^%
\infty[0,b]$. We take the scalar product of both sides of (\ref{rec:eq:kernel-pde}) with each Legendre polynomial and use (\ref{rec:eq:kernel-legendre-expansion}). Approximating $p$ and $q$ by smooth functions and passing to the limit shows that these coefficient formulas remain valid under the stated regularity assumptions.

Set
\begin{equation*}
\mathcal{F}_\pm(x,u) =xK_\pm(x,xu), \qquad-1\leq u\leq1.
\end{equation*}
Then
\begin{equation}
\mathcal{F}_\pm(x,u) =\sum_{n=0}^{\infty}\kappa_{n}^\pm(x)P_{n}(u).
\label{rec:eq:F-Legendre}
\end{equation}
A change of variables in (\ref{rec:eq:kernel-pde}) gives
\begin{equation}
\begin{split}
0={} & x^{2}\frac{\partial^{2}\mathcal{F}_\pm}{\partial x^{2}} -2x\frac{%
\partial \mathcal{F}_\pm}{\partial x} -2xu\frac{\partial^{2}\mathcal{F}_\pm}{%
\partial x\partial u} +(u^{2}-1)\frac{\partial^{2}\mathcal{F}_\pm}{\partial
u^{2}} \\
& +4u\frac{\partial \mathcal{F}_\pm}{\partial u} +2\mathcal{F}_\pm \mp 2 i
p(x)x\frac{\partial \mathcal{F}_\pm}{\partial u} -q(x)x^{2}\mathcal{F}_\pm.
\end{split}
\label{rec:eq:F-pde}
\end{equation}

We use the standard identities
\begin{align}
P_{m}^{\prime}(u) &=\sum_{%
\substack{0\leq k\leq m-1 \\ m-k\text{ odd}}
} (2k+1)P_{k}(u),  \label{rec:eq:Legendre-derivative-expansion}\\
uP_{m}^{\prime}(u) &=mP_{m}(u) +\sum_{%
\substack{ 0\leq k\leq m-2 \\  m-k\text{ even}}
} (2k+1)P_{k}(u).  \label{rec:eq:u-Legendre-derivative-expansion}
\end{align}
These identities hold for $m\geq1$. An empty sum equals zero.

Define
\begin{align*}
A_{n}^{\pm }&=x^{2}\left( \kappa _{n}^{\pm }\right) ^{\prime \prime
}-2x\left( \kappa _{n}^{\pm }\right) ^{\prime }+\left(
n(n+1)+2-q(x)x^{2}\right) \kappa _{n}^{\pm },\\
B_{n}^{\pm }&=\kappa _{n}^{\pm }-x\left( \kappa _{n}^{\pm }\right) ^{\prime }.
\end{align*}

Define the scaled coefficients
\begin{equation}
\widetilde{\kappa}_{n}^\pm(x) =x^{n}\kappa_{n}^\pm(x), \qquad n=0,1,2,\ldots.
\label{rec:eq:scaled-kappa-definition}
\end{equation}

\begin{theorem}
Under the smoothness assumption of this subsection, for every $n\geq2$ the
scaled coefficients satisfy
\begin{equation}
\frac{1}{2n+1}L\left[ \widetilde{\kappa}_{n}^\pm\right] = \frac{x^{2}}{%
2n-3}L\left[ \widetilde{\kappa}_{n-2}^\pm\right]
 -2x\left( \left( \widetilde{\kappa}_{n-2}^\pm\right) ^{\prime}-\frac{n-1}{x%
}\widetilde{\kappa}_{n-2}^\pm \right) \mp 2 i p(x)\widetilde{\kappa}%
_{n-1}^\pm.
\label{rec:eq:scaled-kappa-differential-recurrence}
\end{equation}
For $n\geq1$, the scaled coefficients satisfy the initial conditions
\begin{equation}
\widetilde{\kappa}_{n}^\pm(0)=0, \qquad\left( \widetilde{\kappa}%
_{n}^\pm\right) ^{\prime}(0)=0.  \label{rec:eq:scaled-kappa-zero-initial}
\end{equation}
\end{theorem}

\noindent \textit{Proof.} For an integrable function $v$ on $[-1,1]$, denote
its $n$th Legendre coefficient by
\begin{equation*}
\lbrack v]_{n}=\frac{2n+1}{2}\int_{-1}^{1}v(u)P_{n}(u)\,du.
\end{equation*}%
Thus $[\mathcal{F}_{\pm }]_{n}=\kappa _{n}^{\pm }$. Under the smoothness
assumptions of this subsection, the differentiated Fourier--Legendre partial
sums converge in the required $L^{2}(-1,1)$ spaces. We may therefore pass to
the limit in the scalar products below. Equations (\ref%
{rec:eq:Legendre-derivative-expansion}) and (\ref%
{rec:eq:u-Legendre-derivative-expansion}) imply
\begin{equation*}
\begin{split}
\lbrack \mathcal{F}_{\pm ,u}]_{n}& =(2n+1)\sum_{j\geq 0}\kappa
_{n+1+2j}^{\pm }, \\
\lbrack u\mathcal{F}_{\pm ,u}]_{n}& =n\kappa _{n}^{\pm }+(2n+1)\sum_{j\geq
1}\kappa _{n+2j}^{\pm }.
\end{split}%
\end{equation*}%
The second identity applied to $\mathcal{F}_{\pm ,x}$ replaces each $\kappa
_{m}^{\pm }$ by $(\kappa _{m}^{\pm })^{\prime }$. Furthermore, the Legendre
differential equation implies
\begin{equation*}
(u^{2}-1)P_{m}^{\prime \prime }+4uP_{m}^{\prime
}=m(m+1)P_{m}+2uP_{m}^{\prime }.
\end{equation*}%
Take the scalar product of both sides of (\ref{rec:eq:F-pde}) with $P_{n}$ in $L^{2}(-1,1)$
and multiply the resulting equality by $(2n+1)/2$. Substitution of the
preceding identities and collection of the terms with index $n$ gives
\begin{equation*}
0=A_{n}^{\pm }+2nB_{n}^{\pm }+2(2n+1)\sum_{j\geq 1}B_{n+2j}^{\pm }\mp
2ip(x)x(2n+1)\sum_{j\geq 0}\kappa _{n+1+2j}^{\pm }.
\end{equation*}%
Divide this equality by $2n+1$, and subtract from it the corresponding
equality with index $n-2$, divided by $2n-3$. The terms with indices $%
n+2,n+4,\ldots $ cancel in the sums involving $B_{j}^{\pm }$, while the
terms with indices $n+1,n+3,\ldots $ cancel in the sums involving $\kappa
_{j}^{\pm }$. The remaining equality has the form
\begin{equation*}
0=\frac{A_{n}^{\pm }-2(n+1)B_{n}^{\pm }}{2n+1}-\frac{A_{n-2}^{\pm
}+2(n-2)B_{n-2}^{\pm }}{2n-3}\pm 2ip(x)x\kappa _{n-1}^{\pm }.
\end{equation*}%
Expanding the derivatives in the definitions of $A_{n}^{\pm }$ and $%
B_{n}^{\pm }$, we obtain
\begin{equation*}
A_{n}^{\pm }-2(n+1)B_{n}^{\pm }=x^{2-n}L[x^{n}\kappa _{n}^{\pm }],\qquad
A_{n-2}^{\pm }+2(n-2)B_{n-2}^{\pm }=x^{n+1}L[\kappa _{n-2}^{\pm }/x^{n-1}].
\end{equation*}%
Since $\widetilde{\kappa }_{j}^{\pm }=x^{j}\kappa _{j}^{\pm }$, multiply the
preceding equality by $x^{n-2}$ and use the identity
\begin{equation*}
x^{2n-1}L[\widetilde{\kappa }_{n-2}^{\pm }/x^{2n-3}]=x^{2}L[\widetilde{%
\kappa }_{n-2}^{\pm }]-2(2n-3)x\left( (\widetilde{\kappa }_{n-2}^{\pm
})^{\prime }-\frac{n-1}{x}\widetilde{\kappa }_{n-2}^{\pm }\right) ,
\end{equation*}%
to substitute for the two expressions containing $L$. The result is (\ref%
{rec:eq:scaled-kappa-differential-recurrence}).

It remains to verify the initial conditions. Continuity of $K_\pm$ near $%
(0,0)$ gives, for each fixed $n$,
\begin{equation*}
|\kappa_n^\pm(x)| \leq (2n+1)x\max_{|t|\leq x}|K_\pm(x,t)|=O(x)
\quad(x\to0+).
\end{equation*}
Consequently $\widetilde{\kappa}_n^\pm(x)=O(x^{n+1})$. For $n\ge1$ this
implies (\ref{rec:eq:scaled-kappa-zero-initial}), because the scaled
coefficients are differentiable at zero under the smoothness assumptions of
this subsection. \hfill$\square$

\begin{remark}
The first two terms in (\ref{rec:eq:scaled-kappa-differential-recurrence})
are the Schr\"odinger recurrence of \cite{KNT2017}. The last term couples
the even and odd coefficient subsequences and vanishes when $p\equiv0$.
\end{remark}

\subsection{Recurrent integration procedure for computing NSBF coefficients}

\label{rec:sec:recursive-integration}

Equation (\ref{rec:eq:scaled-kappa-differential-recurrence}) is already a
recursive sequence of inhomogeneous equations for the differential
expression $L$. Direct use of that recurrence, however, would require
differentiating previously computed coefficients, which is inconvenient for
numerical implementation. We therefore eliminate those derivatives using the
recursive integration procedure of \cite[Section~6]{KNT2017}.

Let
\begin{equation*}
f(x)=c(0,x).
\end{equation*}%
Since $e_{+}(0,x)$ and $e_{-}(0,x)$ satisfy the same equation and the same
initial conditions, one has
\begin{equation*}
e_{+}(0,x)=e_{-}(0,x)=f(x),
\end{equation*}%
where
\begin{equation}
L[f]=0,\qquad f(0)=1,\qquad f^{\prime }(0)=0.  \label{rec:eq:f-initial}
\end{equation}%
Let $g=s(0,\cdot )$ be the second solution for $\rho =0$, so that
\begin{equation}
L[g]=0,\qquad g(0)=0,\qquad g^{\prime }(0)=1.  \label{rec:eq:g-initial}
\end{equation}

Let $F$ be a nonvanishing auxiliary solution of
\begin{equation}
L[F]=0  \label{rec:eq:auxiliary-equation}
\end{equation}%
normalized by
\begin{equation*}
F(0)=1.
\end{equation*}%
Such a solution always exists on $[0,b]$, also for complex $q$. Indeed, $%
f,g\in C^{1}[0,b]$, and their Wronskian equals one by (\ref{rec:eq:f-initial}%
)--(\ref{rec:eq:g-initial}), so they have no common zero. For real $q$ one
can take $F=f+ig$. For complex $q$, choose some $\omega \in \mathbb{C}$ outside the set $%
\{-f(x)/g(x):g(x)\neq 0,\ 0\leq x\leq b\}$ and take $F=f+\omega g$. Then $F$
has no zeros and $F(0)=1$. The existence of a nonvanishing linear
combination for complex coefficients was proved in \cite[Remark~5]%
{KravchenkoPorter2010}, see also \cite[Corollary~2.3]{CamporesiDiScala2011}.

The choice of $F$ does not affect the resulting coefficients.

Equation (\ref{rec:eq:auxiliary-equation}) implies the P\'olya factorization
\begin{equation}
L[y]=\frac{1}{F}\frac{d}{dx} \left( F^{2}\frac{d}{dx}\left( \frac{y}{F}%
\right) \right).  \label{rec:eq:Polya-factorization}
\end{equation}
This factorization turns the differential recurrence into successive
integrations.

For each choice of sign, suppose that $\widetilde{\kappa}_{0}^\pm$ and $%
\widetilde{\kappa}_{1}^\pm$ have already been computed. For $n\geq2$, define
\begin{align}
\eta_{n}^\pm(x) &=\int_{0}^{x}\left(tF^{\prime }(t)+(n-1)F(t)\right)
\widetilde{\kappa}_{n-2}^\pm(t)\,dt,  \label{rec:eq:eta-n}\\
\zeta_{n}^\pm(x) &=\int_{0}^{x}F(t)p(t)\widetilde{\kappa}_{n-1}^\pm(t)\,dt,
\label{rec:eq:zeta-n}\\
\theta_{n}^\pm(x) &=\int_{0}^{x}\frac{\eta_{n}^\pm(t) -tF(t)\widetilde{\kappa}%
_{n-2}^\pm(t)}{F(t)^{2}}\,dt,  \label{rec:eq:theta-n}\\
\mu_{n}^\pm(x) &=\int_{0}^{x}\frac{\zeta_{n}^\pm(t)}{F(t)^{2}}\,dt.
\label{rec:eq:mu-n}
\end{align}
Differentiating the definitions (\ref{rec:eq:eta-n})--(\ref{rec:eq:mu-n})
gives
\begin{align}
\left( \theta_{n}^\pm\right) ^{\prime}(x) &=\frac{\eta _{n}^\pm(x)-xF(x)%
\widetilde{\kappa}_{n-2}^\pm(x)}{F(x)^{2}}.  \label{rec:eq:theta-prime}\\
\left( \mu_{n}^\pm\right) ^{\prime}(x) &=\frac{\zeta _{n}^\pm(x)}{F(x)^{2}}.
\label{rec:eq:mu-prime}
\end{align}

\begin{theorem}
Assume $q\in L^2(0,b)$ and $p\in W^{1,2}(0,b)$, and let $F$ be a solution of
(\ref{rec:eq:auxiliary-equation}) without zeros on $[0,b]$, normalized by $%
F(0)=1$. For $n=2,3,\ldots$ and $0<x\leq b$, the scaled coefficients can be
obtained from the following recursive integration formula, which requires no differentiation:
\begin{equation}
\widetilde{\kappa}_{n}^\pm(x) = \frac {2n+1}{2n-3} \left( x^{2}%
\widetilde{\kappa}_{n-2}^\pm(x) +2(2n-1)F(x)\theta _{n}^\pm(x) \right)
\mp 2 i(2n+1)F(x)\mu_{n}^\pm(x).
\label{rec:eq:main-recursive-integration}
\end{equation}
The NSBF coefficients are then
\begin{equation*}
\kappa_{n}^\pm(x)=\frac{\widetilde{\kappa}_{n}^\pm(x)}{x^{n}}, \qquad x>0.
\end{equation*}
\end{theorem}

\noindent \textit{Proof.} Set $w=\widetilde{\kappa }_{n-2}^{\pm }$.
Equations (\ref{rec:eq:eta-n})--(\ref{rec:eq:mu-n}) and the P\'{o}lya
factorization (\ref{rec:eq:Polya-factorization}) give
\begin{equation*}
(F^{2}(\theta _{n}^{\pm })^{\prime })^{\prime }=F((n-2)w-xw^{\prime
}),\qquad (F^{2}(\mu _{n}^{\pm })^{\prime })^{\prime }=Fp\widetilde{\kappa }%
_{n-1}^{\pm },
\end{equation*}%
and hence
\begin{equation*}
L[F\theta _{n}^{\pm }]=(n-2)w-xw^{\prime },\qquad L[F\mu _{n}^{\pm }]=p%
\widetilde{\kappa }_{n-1}^{\pm }.
\end{equation*}%
Together with $L[x^{2}w]=x^{2}L[w]+4xw^{\prime }+2w$, these identities give
the equality
\begin{equation*}
L[\widetilde{\kappa }_{n}^{\pm }]= \frac{2n+1}{2n-3}\left\{
x^{2}L[w]-2(2n-3)x\left( w^{\prime }-\frac{n-1}{x}w\right) \right\}
 \mp 2i(2n+1)p\widetilde{\kappa }_{n-1}^{\pm }.
\end{equation*}%
After division by $2n+1$, this is precisely (\ref%
{rec:eq:scaled-kappa-differential-recurrence}). All four integrals in (\ref%
{rec:eq:eta-n})--(\ref{rec:eq:mu-n}) vanish at the origin. Together with the
initial conditions of the previously computed scaled coefficients, this
shows that the right-hand side of (\ref{rec:eq:main-recursive-integration})
and its first derivative vanish at $x=0$. Uniqueness for the initial value
problem with these two zero initial values proves the formula.

For general $q\in L^2(0,b)$ and $p\in W^{1,2}(0,b)$, take smooth $q_\nu\to q$
in $L^2$ and $p_\nu\to p$ in $W^{1,2}$. Let $F_\nu$ solve $F_\nu^{\prime
\prime }-q_\nu F_\nu=0$ with the same initial values as $F$. On a closed
interval on which $F$ has no zero, $\min|F|>0$. The approximating solutions
and their first derivatives converge uniformly, so the denominators $F_\nu^2$
are bounded away from zero for all sufficiently large approximation indices $%
\nu$. Uniform convergence of the corresponding kernels gives convergence of
the Fourier--Legendre coefficients. Induction on the coefficient index now
permits passage to the limit in each of the four successive integrals and in
(\ref{rec:eq:main-recursive-integration}). \hfill$\square$

\begin{remark}
Each step uses four successive integrals and algebraic operations at each
grid point, and requires no differentiation. When $p=0$, the two $p$%
-dependent quantities $\zeta _{n}^{\pm }$ and $\mu _{n}^{\pm }$ vanish, and
we obtain the recursive integration procedure from \cite{KNT2017}.
\end{remark}

\subsection{The initial coefficients}

\label{rec:sec:initial-coefficients}

The recurrence (\ref{rec:eq:main-recursive-integration}) starts from the
coefficients of indices zero and one. We now express these initial
coefficients in terms of solutions for $\rho =0$ and successive integrals.
They also provide the data for the reconstruction of $p$ and $q$, which is
especially convenient for solving corresponding inverse problems.

Recall $f$ and $g$ from (\ref{rec:eq:f-initial})--(\ref{rec:eq:g-initial}).
Let $c_{1}=\left. \partial _{\rho }c(\rho ,\cdot )\right\vert _{\rho =0}$.
Differentiating the pencil with respect to $\rho $ and then setting $\rho =0$
gives
\begin{equation}
L[c_{1}]=2pf,\qquad c_{1}(0)=0,\qquad c_{1}^{\prime }(0)=0.
\label{rec:eq:c1-equation}
\end{equation}%
Define
\begin{equation*}
A(x)=\int_{0}^{x}p(t)f(t)^{2}\,dt,\qquad B(x)=\int_{0}^{x}p(t)f(t)g(t)\,dt.
\end{equation*}%
Then
\begin{equation}
c_{1}(x)=2\left( g(x)A(x)-f(x)B(x)\right) ,\qquad c_{1}^{\prime }(x)=2\left(
g^{\prime }(x)A(x)-f^{\prime }(x)B(x)\right) .
\label{rec:eq:c1-global-quadrature}
\end{equation}%
The terms containing $A^{\prime }$ and $B^{\prime }$ cancel in the
derivative formula.

\begin{proposition}
The first coefficients in (\ref{rec:eq:e-pm-nsbf}) have the form
\begin{align}
\kappa _{0}^{+}&=\tfrac{1}{2}(f-e^{-i\sigma }),  & \kappa_{0}^{-}&=\tfrac12(f-e^{i\sigma}), \label{rec:eq:alpha0beta0}\\
\kappa_{1}^{+}&=\tfrac32\left((g-ic_1)/x-e^{-i\sigma}\right), & \kappa _{1}^{-}&=\tfrac{3}{2}\left( (g+ic_{1})/x-e^{i\sigma }\right).\label{rec:eq:alpha1beta1}
\end{align}
Equivalently, the scaled initial coefficients are
\begin{align}
\widetilde{\kappa }_{0}^{\pm }&=\tfrac{1}{2}(f-e^{\mp i\sigma }),
\label{rec:eq:scaled-kappa0-initial}\\
\widetilde{\kappa }_{1}^{\pm }&=\tfrac{3}{2}(g\mp ic_{1}-xe^{\mp i\sigma }).
\label{rec:eq:scaled-kappa1-initial}
\end{align}%
The first coefficients for $c$ and $s$ have the form
\begin{align}
a_{0}&=f-\cos \sigma , & u_0&=\sin\sigma, \label{rec:eq:a0u0}\\
v_0&=3(g/x-\cos\sigma),& b_0&=3(c_1/x-\sin\sigma). \label{rec:eq:v0b0}
\end{align}
\end{proposition}

\noindent\textit{Proof.} Since $P_{0}(z)=1$, formula (\ref%
{rec:eq:direct-kappa-formula}) gives
\begin{equation*}
\kappa_{0}^\pm =\frac{1}{2}\left( \Phi_{0}^\pm-e^{\mp i\sigma} \right) .
\end{equation*}
Equation (\ref{rec:eq:Phi-recurrence}) and the initial data show that $%
\Phi_{0}^\pm=f$, which proves (\ref{rec:eq:alpha0beta0}).

For $k=1$, equation (\ref{rec:eq:Phi-recurrence}) is
\begin{equation*}
L\left[ \Phi_{1}^\pm\right] =\mp 2 i p f, \qquad\Phi _{1}^\pm(0)=0,
\qquad\left( \Phi_{1}^\pm\right) ^{\prime}(0)=1.
\end{equation*}
Equations (\ref{rec:eq:g-initial}) and (\ref{rec:eq:c1-equation}) imply
\begin{equation*}
\Phi_{1}^\pm=g\mp i c_1.
\end{equation*}
Since $P_{1}(z)=z$, formula (\ref{rec:eq:direct-kappa-formula}) gives
\begin{equation*}
\kappa_{1}^\pm =\frac{3}{2} \left( \frac{\Phi_{1}^\pm}{x} -e^{\mp i\sigma}
\right) .
\end{equation*}
Equations (\ref{rec:eq:alpha1beta1}) follow. Taking
their sums and differences gives (\ref{rec:eq:a0u0})--(\ref%
{rec:eq:v0b0}). \hfill$\square$

\begin{remark}
The recurrence can be started at $n=1$ by introducing the auxiliary
functions
\[
\widetilde{\kappa}_{-1}^{\pm}(x)
=-\frac{\kappa_0^\pm(x)}{x},\qquad x>0.
\]
Since $\kappa_0^\pm(0)=0$ and
$(\kappa_0^\pm)'(0)=\pm ip(0)/2$, these functions extend continuously
to zero with
\[
\widetilde{\kappa}_{-1}^{\pm}(0)=\mp\frac{i}{2}p(0).
\]
With this definition, both
(\ref{rec:eq:scaled-kappa-differential-recurrence}) and
(\ref{rec:eq:main-recursive-integration}) remain valid for $n=1$.
Indeed, substitution verifies the differential recurrence, and the
integration formula follows from the same zero initial conditions.
Thus only $\kappa_0^\pm$ requires separate initialization.
The negative index is used solely in the recurrence; the
Fourier--Legendre and NSBF expansions still start at $n=0$.
\end{remark}

\section{Recovery of $p$ and $q$ from the first coefficients}

\label{sec:potential-recovery}

In this section we develop two procedures for the recovery of $p$ and $q$
from the first coefficients. The procedures are mutually complementary.

\begin{theorem}[Recovery from three initial coefficients]
\label{thm:initial-coefficient-reconstruction} Assume $q\in L^{2}(0,b)$ and $%
p\in W^{1,2}(0,b)$. Let $I$ be a subinterval of $[0,b]$ containing zero on
which $1-u_{0}(x)^{2}\neq 0$. Denote by $C$ the continuous square root of $%
1-u_{0}^{2}$ on $I$ satisfying $C(0)=1$. The zero energy ($\rho =0$)
solutions can be recovered on $I$ from the first coefficients as follows
\begin{equation*}
f(x)=a_{0}(x)+C(x),\qquad g(x)=x\left( \frac{v_{0}(x)}{3}+C(x)\right) .
\end{equation*}%
They determine the coefficients of (\ref{eq:pencil}) through the formulas
\begin{equation}
\begin{split}
p(x)& =\frac{u_{0}^{\prime }(x)}{C(x)}, \\
q(x)& =f^{\prime \prime }(x)g^{\prime }(x)-f^{\prime }(x)g^{\prime \prime
}(x) \\
&=\frac{\overline{f(x)}f^{\prime \prime }(x)+\overline{g(x)}g^{\prime
\prime }(x)}{|f(x)|^{2}+|g(x)|^{2}}.
\end{split}
\label{rec:eq:recovery-formulas}
\end{equation}%
The formula for $p$ holds at every point of $I$, and those for $q$ hold
almost everywhere on $I$.
\end{theorem}

\noindent \textit{Proof.} By (\ref{rec:eq:a0u0}), $u_{0}=\sin \sigma $.
Differentiation of $u_{0}$ proves the formula for $p$. Equations (\ref%
{rec:eq:a0u0}) and (\ref{rec:eq:v0b0}) now give the stated expressions for $f$
and $g$.

The solutions $f$ and $g$ belong to $W^{2,2}(0,b)$ and satisfy
\begin{equation}
Ly=0,  \label{Ly=0}
\end{equation}
almost everywhere. Their Wronskian $fg^{\prime }-f^{\prime }g$ equals one by
(\ref{rec:eq:f-initial})--(\ref{rec:eq:g-initial}). This proves the first
expression for $q$. Moreover, $f$ and $g$ cannot vanish at the same point.
Multiplying their differential equations by $\overline{f}$ and $\overline{g}$%
, respectively, and adding proves the second expression. \hfill $\square $

The two formulas for $q$ in (\ref{rec:eq:recovery-formulas}) remain valid
when either zero energy solution vanishes. The denominator $|f|^{2}+|g|^{2}$
has a positive lower bound on every compact subinterval of $I$. If $q$ is
continuous, both formulas for it hold pointwise.

The formula for $p$ uses the branch $C=\cos \sigma $ fixed at zero. It
applies until a zero of $\cos \sigma $ is reached. If this branch is known
beyond such a zero, the identity $p=Cu_{0}^{\prime }-u_{0}C^{\prime }$
avoids division by $C$.

The following formulas use all four coefficients in (\ref{rec:eq:a0u0})--(\ref%
{rec:eq:v0b0}) and require no choice of a square root. They can be applied at
zeros of $\cos \sigma $ where the preceding theorem does not apply.

\begin{proposition}[Recovery from four initial coefficients]
\label{rec:prop:four-coefficient-recovery} Assume $q\in L^{2}(0,b)$ and $%
p\in W^{1,2}(0,b)$. From the four initial coefficients form the auxiliary
functions
\begin{equation*}
H(x)=x\left( \frac{v_{0}(x)}{3}-a_{0}(x)\right) ,\qquad c_{1}(x)=x\left(
u_{0}(x)+\frac{b_{0}(x)}{3}\right) .
\end{equation*}%
Define
\begin{equation*}
U=H^{\prime \prime }+2a_{0}^{\prime },\qquad V=c_{1}^{\prime \prime
}-2u_{0}^{\prime },\qquad W=a_{0}H-u_{0}c_{1}.
\end{equation*}%
The functions $U$ and $V$ belong to $L^{2}(0,b)$, and $W$ is continuous on $%
[0,b]$. On the open set where $W(x)\neq 0$, we have
\begin{equation}
\begin{split}
q(x)& =\frac{a_{0}(x)U(x)-u_{0}(x)V(x)}{W(x)}, \\
p(x)& =\frac{H(x)V(x)-c_{1}(x)U(x)}{2W(x)}
\end{split}
\label{rec:eq:four-coefficient-recovery}
\end{equation}%
almost everywhere. These formulas remain valid at zeros of $\cos \sigma $
where $W$ does not vanish.
\end{proposition}

\noindent \textit{Proof.} Equations (\ref{rec:eq:a0u0}) and (\ref{rec:eq:v0b0})
give $H=g-xf$ as well as the second expression for $c_{1}$. In particular, $H,c_{1}\in W^{2,2}(0,b)$%
. Using the fact that $f$ and $g$ solve (\ref{Ly=0}), and $a_{0}^{\prime
}=f^{\prime }+pu_{0}$, we obtain
\begin{equation*}
U=H^{\prime \prime }+2a_{0}^{\prime }=qH+2pu_{0}.
\end{equation*}%
Similarly, (\ref{rec:eq:c1-equation}) and $u_{0}^{\prime }=p\cos \sigma $
yield
\begin{equation*}
V=c_{1}^{\prime \prime }-2u_{0}^{\prime }=qc_{1}+2pa_{0}.
\end{equation*}%
Both identities hold almost everywhere. This is a linear system for $q$ and $%
p$ with determinant $2(a_{0}H-u_{0}c_{1})=2W$. Solving it gives (\ref%
{rec:eq:four-coefficient-recovery}). \hfill $\square $

One always has $W(0)=0$.
When $p\equiv 0$, the preceding theorem recovers $q$ throughout $[0,b]$,
without any condition on $W$. On the other hand, if $p(0)\neq 0$, then
\begin{equation*}
W(x)=-p(0)^{2}x^{3}+o(x^{3})\qquad (x\downarrow 0).
\end{equation*}%
Indeed, $u_{0}(x)=p(0)x+o(x)$ and $c_{1}(x)=p(0)x^{2}+o(x^{2})$, while $%
a_{0}(x)H(x)=o(x^{3})$ by the initial conditions. In this case $W(x)\neq 0$
for all sufficiently small $x>0$. For example, when $q=0$ and $p=1$, one has
$H=0$, $c_{1}=x^{2}$, and $W=-x^{2}\sin x$. Thus (\ref%
{rec:eq:four-coefficient-recovery}) also applies at $x=\pi /2$, provided $%
b>\pi /2$, although $\cos \sigma $ vanishes there.

Nevertheless, it should be noted that the assumptions do not exclude simultaneous
zeros of $W$ and $\cos\sigma$. Consequently, the two recovery
procedures do not necessarily cover the whole interval.

\section{NSBF representations for derivatives}

\label{der:sec:derivative-nsbf}

To obtain NSBF representations for the first derivatives of the solutions,
we expand the partial derivative of the kernel with respect to its first
variable. Here and below, this derivative is taken at fixed $t$. By Lemma~%
\ref{lem:goursat-sobolev-regularity} with $m=0$, the hypotheses of Theorem~%
\ref{thm:transmutation-representation} imply $K_{\pm ,x}(x,\cdot )\in
L^{2}(-x,x)$ for every $x>0$. We therefore write its Fourier--Legendre
expansion as
\begin{equation}
K_{\pm ,x}(x,t)=\sum_{n=0}^{\infty }\frac{\chi _{n}^{\pm }(x)}{x}P_{n}\left(
\frac{t}{x}\right) ,  \label{der:eq:kernel-x-Legendre}
\end{equation}%
where
\begin{equation}
\chi _{n}^{\pm }(x)=\frac{2n+1}{2}\int_{-x}^{x}K_{\pm ,x}(x,t)P_{n}(t/x)\,dt.
\label{der:eq:chi-definition}
\end{equation}

By Lemma \ref{lem:goursat-sobolev-regularity}
\[
|\chi_n^\pm(x)|\leq \left(\frac{(2n+1)x}{2}\right)^{1/2} \|K_{\pm,x}(x,\cdot)\|_{L^2(-x,x)}\longrightarrow0,\quad x\to 0,
\]
hence the coefficients $\chi_n^\pm$ can be continuously extended to 0 with the limit values
\[
\chi_n^\pm(0)=0,\qquad n=0,1,2,\ldots.
\]

\begin{theorem}
For every $0<x\leq b$ and $\rho \in \mathbb{C}$,
\begin{equation}
e_{\pm }^{\prime }(\rho ,x)= \left( \pm i\rho +\frac{\mathcal{R}(x)}{2}%
\mp \frac{i}{2}p(x)\right) e^{\pm i\vartheta (\rho ,x)}  \pm \frac{i}{2}p(0)e^{\mp i\vartheta (\rho ,x)}+2\sum_{n=0}^{\infty }(\pm
i)^{n}\chi _{n}^{\pm }(x)j_{n}(\rho x).
\label{der:eq:e-derivative-nsbf}
\end{equation}%
At $x=0$, the formula gives $e_\pm'(\rho,0)=\pm i\rho$. For fixed $x$, the
series converge absolutely and uniformly with respect to $\rho $ on every
horizontal strip.
\end{theorem}

\noindent \textit{Proof.} Differentiation of (\ref%
{rec:eq:transmutation-representation}) gives the integral representation

\begin{equation}
\begin{split}
e_{\pm }^{\prime }(\rho ,x)={} & \pm i\left( \rho -p(x)\right) e^{\pm
i\vartheta (\rho ,x)}+K_{\pm }(x,x)e^{\pm i\rho x} \\
& +K_{\pm }(x,-x)e^{\mp i\rho x}+\int_{-x}^{x}K_{\pm ,x}(x,t)e^{\pm i\rho
t}\,dt.%
\end{split}
\label{der:eq:e-derivative-before-characteristics}
\end{equation}%
Using (\ref{rec:eq:kernel-upper-value}), (\ref{rec:eq:kernel-lower-value})
and identity (\ref{eq:legendre-bessel-integral}) leads to (\ref%
{der:eq:e-derivative-nsbf}). The passage from Fourier--Legendre partial sums
to the integral follows from their $L^{2}(-x,x)$ convergence and the
Cauchy--Schwarz inequality. Applying the same inequality to the
Fourier--Legendre coefficients of $K_{\pm ,x}(x,\cdot )$ and using identity (%
\ref{err:eq:bessel-parseval}) proves absolute and uniform convergence on
each horizontal strip. \hfill $\square $

Denote the derivative coefficient families by
\begin{equation}
\widehat{a}_{n}=\chi _{2n}^{+}+\chi _{2n}^{-},\qquad \widehat{b}_{n}=i\left(
\chi _{2n+1}^{+}-\chi _{2n+1}^{-}\right) ,  \label{der:eq:ahat-bhat}
\end{equation}%
\begin{equation}
\widehat{u}_{n}=-i\left( \chi _{2n}^{+}-\chi _{2n}^{-}\right) ,\qquad
\widehat{v}_{n}=\chi _{2n+1}^{+}+\chi _{2n+1}^{-}.  \label{der:eq:uhat-vhat}
\end{equation}

\begin{theorem}
For every $0<x\leq b$ and $\rho \in \mathbb{C}$, the derivatives of the
fundamental solutions admit the representations
\begin{equation}
\begin{split}
c^{\prime }(\rho ,x)={}& -\rho \sin \vartheta (\rho ,x)+\frac{\mathcal{R}(x)%
}{2}\cos \vartheta (\rho ,x)+\frac{p(x)+p(0)}{2}\sin \vartheta (\rho ,x) \\
& +\sum_{n=0}^{\infty }(-1)^{n}\widehat{a}_{n}(x)j_{2n}(\rho
x)+\sum_{n=0}^{\infty }(-1)^{n}\widehat{b}_{n}(x)j_{2n+1}(\rho x),
\end{split}
\label{der:eq:c-derivative-nsbf}
\end{equation}%
\begin{equation}
\begin{split}
s^{\prime }(\rho ,x)={}& \cos \vartheta (\rho ,x) \\
& +\frac{\mathcal{R}(x)\sin \vartheta (\rho ,x)+(p(0)-p(x))\cos \vartheta
(\rho ,x)+2\widehat{u}_{0}(x)j_{0}(\rho x)}{2\rho } \\
& +\frac{1}{\rho }\sum_{n=1}^{\infty }(-1)^{n}\widehat{u}_{n}(x)j_{2n}(\rho
x)+\frac{1}{\rho }\sum_{n=0}^{\infty }(-1)^{n}\widehat{v}_{n}(x)j_{2n+1}(%
\rho x).
\end{split}
\label{der:eq:s-derivative-nsbf}
\end{equation}%
For each fixed $x>0$, the two series in each representation converge
absolutely and uniformly with respect to $\rho $ on every horizontal strip $|%
\operatorname{Im}\rho |\leq C$, $C\geq 0$. In (\ref{der:eq:s-derivative-nsbf}), both
the finite quotient and each divided Bessel function term $j_{k}(\rho
x)/\rho $ are assigned their corresponding asymptotic values at $\rho =0$.
\end{theorem}

\noindent\textit{Proof.} For $\rho\neq0$, the representations follow from (%
\ref{der:eq:e-derivative-nsbf}) and the coefficient definitions (\ref%
{der:eq:ahat-bhat})--(\ref{der:eq:uhat-vhat}). We have placed the term
containing $\widehat u_0$ together with the explicit terms divided by $\rho$
in order to account for their cancellation at zero.

Fix $x>0$. Parseval's identity for (\ref{der:eq:kernel-x-Legendre}) gives
\begin{equation*}
\sum_{n=0}^{\infty }\frac{|\chi _{n}^{\pm }(x)|^{2}}{2n+1}=\frac{x}{2}\Vert
K_{\pm ,x}(x,\cdot )\Vert _{L^{2}(-x,x)}^{2}<\infty .
\end{equation*}
The absolute and uniform convergence in (\ref{der:eq:c-derivative-nsbf})
follows from the preceding theorem and the coefficient definitions. For the
divided series in (\ref{der:eq:s-derivative-nsbf}), consider
\begin{equation*}
\psi _{\rho }(t)=\frac{e^{i\rho t}-1}{\rho }=it\int_{0}^{1}e^{i\theta \rho
t}\,d\theta ,\qquad \psi _{0}(t)=it.
\end{equation*}%
On $|\operatorname{Im}\rho |\leq C$, we have $|\psi _{\rho }(t)|\leq |t|e^{C|t|}$,
and hence $\Vert \psi _{\rho }\Vert _{L^{2}(-x,x)}\leq \Theta _{C}(x)$ by (%
\ref{err:eq:Theta}). For $n\geq 1$, orthogonality of $P_{n}$ to constants
and (\ref{eq:legendre-bessel-integral}) imply
\begin{equation*}
\int_{-x}^{x}\psi _{\rho }(t)P_{n}(t/x)\,dt=2xi^{n}\frac{j_{n}(\rho x)}{\rho
}.
\end{equation*}%
Bessel's inequality therefore yields
\begin{equation*}
\sum_{n=1}^{\infty }(2n+1)\left\vert \frac{j_{n}(\rho x)}{\rho }\right\vert
^{2}\leq \frac{\Theta _{C}(x)^{2}}{2x}.
\end{equation*}%
For $N\geq 0$, the Cauchy--Schwarz inequality now gives
\begin{equation*}
\sum_{n>N}\left\vert \chi _{n}^{\pm }(x)\frac{j_{n}(\rho x)}{\rho }%
\right\vert \leq \frac{\Theta _{C}(x)}{\sqrt{2x}}\left( \sum_{n>N}\frac{%
|\chi _{n}^{\pm }(x)|^{2}}{2n+1}\right) ^{1/2}.
\end{equation*}%
The right-hand side tends to zero independently of $\rho $ in the strip. By (%
\ref{der:eq:uhat-vhat}), this bounds the absolute tails of both divided
series in (\ref{der:eq:s-derivative-nsbf}) and proves their absolute and
uniform convergence.

Finally, evaluate the integral representation (\ref%
{der:eq:e-derivative-before-characteristics}) at $\rho =0$ and use (\ref%
{rec:eq:kernel-upper-value})--(\ref{rec:eq:kernel-lower-value}). Since $%
e_{+}^{\prime }(0,x)=e_{-}^{\prime }(0,x)$ and the kernel integrals equal $%
2\chi _{0}^{\pm }(x)$, subtraction gives
\begin{equation*}
2\widehat{u}_{0}(x)=\mathcal{R}(x)\sin \sigma (x)+(p(x)-p(0))\cos \sigma (x).
\end{equation*}%
Thus the numerator of the finite quotient in (\ref{der:eq:s-derivative-nsbf}%
) vanishes at $\rho =0$. This quotient has a removable singularity there.
The divided Bessel function terms also extend through zero, with $j_{1}(\rho
x)/\rho \rightarrow x/3$ and $j_{n}(\rho x)/\rho \rightarrow 0$ for $n\geq 2$%
. The uniform convergence just proved permits evaluation term by term at
zero. \hfill $\square $

\bigskip

Note that for $\rho =0$, from (\ref{der:eq:c-derivative-nsbf}) and (\ref%
{der:eq:s-derivative-nsbf}) we obtain

\begin{equation*}
c^{\prime }(0,x)=\frac{\mathcal{R}(x)}{2}\cos \sigma (x)-\frac{p(x)+p(0)}{2}%
\sin \sigma (x)+\widehat{a}_{0}(x),
\end{equation*}

\begin{equation*}
s^{\prime }(0,x)=\cos \sigma (x)+\frac{x}{2}\left( \mathcal{R}(x)\cos \sigma
(x)+(p(0)-p(x))\sin \sigma (x)\right) +\frac{x}{3}\widehat{v}_{0}(x).
\end{equation*}

\subsection{\protect\bigskip Convergence of the representations for
derivatives}
\label{der:sec:convergence}
For $N\geq 0$, let $K_{\pm,x,N}$ be the $N$th Fourier--Legendre partial sum of $K_{\pm,x}$,
\begin{equation*}
K_{\pm ,x,N}(x,t)=\sum_{n=0}^{N}\frac{\chi _{n}^{\pm }(x)}{x}%
P_{n}(t/x),
\end{equation*}
and denote the $L^2$ norm of the remainder by $\widehat d_{\pm,N}$:
\begin{equation*}
\widehat{d}_{\pm ,N}(x)=\Vert K_{\pm ,x}(x,\cdot )-K_{\pm
,x,N}(x,\cdot )\Vert _{L^{2}(-x,x)}.
\end{equation*}%
Let $e_{\pm ,N}^{[1]}$ denote the expression obtained from (\ref%
{der:eq:e-derivative-nsbf}) by retaining its explicit terms and the Bessel
function terms with $0\leq n\leq N$. Define
\begin{equation*}
c_{N}^{[1]}=\frac{e_{+,N}^{[1]}+e_{-,N}^{[1]}}{2},\qquad s_{N}^{[1]}=\frac{%
e_{+,N}^{[1]}-e_{-,N}^{[1]}}{2i\rho },
\end{equation*}%
At $\rho =0$, we have
\begin{equation*}
e_{\pm ,N}^{[1]}(0,x)=e_{\pm }^{\prime }(0,x)=c_{N}^{[1]}(0,x)=c^{\prime
}(0,x),\qquad s_{N}^{[1]}(0,x)=%
\begin{cases}
s^{\prime }(0,x)-x\widehat{v}_{0}(x)/3, & N=0, \\
s^{\prime }(0,x), & N\geq 1.%
\end{cases}%
\end{equation*}%
Thus $s_{N}^{[1]}$ is defined at zero without dividing by $\rho $.

\begin{theorem}[Derivative NSBF truncation estimates]
\label{thm:derivative-error-estimates} Assume $q\in L^{2}(0,b)$ and $p\in
W^{1,2}(0,b)$. Let $N\geq 0$, $%
0<x\leq b$, $C\geq 0$, and $|\operatorname{Im}\rho |\leq C$. Then
\begin{align}
|e_{\pm }^{\prime }(\rho ,x)-e_{\pm ,N}^{[1]}(\rho ,x)|& \leq \Lambda _{C}(x)%
\widehat{d}_{\pm ,N}(x),  \label{der:err:e-prime}\\
|c^{\prime }(\rho ,x)-c_{N}^{[1]}(\rho ,x)|&\leq \frac{\Lambda _{C}(x)}{2}%
\left( \widehat{d}_{+,N}(x)+\widehat{d}_{-,N}(x)\right) ,
\label{der:err:c-prime}
\end{align}
and
\begin{equation}
|\rho |\,|s^{\prime }(\rho ,x)-s_{N}^{[1]}(\rho ,x)|\leq \frac{\Lambda
_{C}(x)}{2}\left( \widehat{d}_{+,N}(x)+\widehat{d}_{-,N}(x)\right) .
\label{der:err:rho-s-prime}
\end{equation}%
For $s^{\prime }$ the following estimate is also valid
\begin{equation}
|s^{\prime }(\rho ,x)-s_{N}^{[1]}(\rho ,x)|\leq \frac{\Theta _{C}(x)}{2}%
\left( \widehat{d}_{+,N}(x)+\widehat{d}_{-,N}(x)\right) .
\label{der:err:s-prime}
\end{equation}%
If, for an integer $r\geq 1$, $q\in W^{r,2}(0,b)$, $p\in W^{r+1,2}(0,b)$, then
\begin{equation*}
\begin{split}
\widehat{\mathcal{A}}_{r}^{\pm }& :=\sup_{0<x\leq b}\Vert \partial
_{t}^{r}K_{\pm ,x}(x,\cdot )\Vert _{L^{2}(-x,x)}, \\
\widehat{\mathcal{A}}_{r}& :=\max \{\widehat{\mathcal{A}}_{r}^{+},\widehat{%
\mathcal{A}}_{r}^{-}\}\leq \widehat{\mathcal{C}}_{r}\!\left( b,\Vert q\Vert
_{W^{r,2}},\Vert p\Vert _{W^{r+1,2}}\right) <\infty .
\end{split}%
\end{equation*}%
Here $\widehat{\mathcal{C}}_{r}$ is locally bounded in its arguments, and
the constants $C_{r}$ below depend only on $r$. For $N>r$,
\begin{equation}
\widehat{d}_{\pm ,N}(x)\leq \frac{C_{r}x^{r}}{N^{r}}\left\Vert \left( 1-%
\frac{t^{2}}{x^{2}}\right) ^{r/2}\partial _{t}^{r}K_{\pm ,x}(x,t)\right\Vert
_{L^{2}(-x,x)}  \label{der:err:kernel-rate}
\end{equation}

\begin{equation}
\widehat d_{\pm,N}(x)\leq\frac{C_r\widehat{\mathcal{A}}_r x^r}{N^r}.
\label{der:err:coefficient-rate}
\end{equation}
Consequently, under these coefficient assumptions, the left-hand sides of (%
\ref{der:err:e-prime}), (\ref{der:err:c-prime}), and (\ref%
{der:err:rho-s-prime}) are bounded by $C_r\widehat{\mathcal{A}}_r
x^r\Lambda_C(x)/N^r$, while that of (\ref{der:err:s-prime}) is bounded by $%
C_r\widehat{\mathcal{A}}_r x^r\Theta_C(x)/N^r$. For every $R>0$ and $%
|\rho|\leq R$, the factorial bounds (\ref{err:eq:e-factorial})--(\ref%
{err:eq:e-over-rho-factorial}) also hold with $e_\pm-e_{\pm,N}$ and $%
d_{\pm,N}$ replaced by $e_\pm^{\prime }-e_{\pm,N}^{[1]}$ and $\widehat
d_{\pm,N}$, respectively.
\end{theorem}

\noindent \textit{Proof.} For $\widehat{R}_{\pm ,N}=K_{\pm ,x}-K_{\pm ,x,N}$%
,
\begin{equation*}
e_{\pm }^{\prime }(\rho ,x)-e_{\pm ,N}^{[1]}(\rho ,x)=\int_{-x}^{x}\widehat{R%
}_{\pm ,N}(x,t)e^{\pm i\rho t}\,dt.
\end{equation*}%
This is the same remainder identity as for $e_{\pm }$, with $K_{\pm }$ and $%
d_{\pm ,N}$ replaced by $K_{\pm ,x}$ and $\widehat{d}_{\pm ,N}$. In
particular, $\widehat{R}_{\pm ,N}$ is orthogonal to $1$. The Cauchy--Schwarz
proof of (\ref{err:eq:e-estimate})--(\ref{err:eq:s-estimate}), with $1$
subtracted before division by $\rho $, proves (\ref{der:err:e-prime})--(\ref%
{der:err:s-prime}).

Lemma~\ref{lem:goursat-sobolev-regularity}, applied with $m=r$, gives the
indicated bound on $\partial _{t}^{r}K_{\pm ,x}$. Applying the same estimate as in Theorem \ref{thm:kernel-projection-rates} for Fourier--Legendre remainders, after rescaling, proves (\ref{der:err:kernel-rate}) and (\ref{der:err:coefficient-rate}).
Finally, $\widehat{R}_{\pm ,N}$ is orthogonal to every polynomial of degree
at most $N$. The Taylor remainder argument in Proposition~\ref%
{prop:compact-factorial-convergence} applies unchanged. \hfill $\square $

\subsection{Direct formulas for the derivative NSBF coefficients}

We use the transmutation operators and transmuted monomials introduced in
Subsection~\ref{rec:sec:transmuted}, and the functions $f,g,c_1$ introduced
in Subsection~\ref{rec:sec:initial-coefficients}. The characteristic values $%
K_\pm(x,x)$ and $K_\pm(x,-x)$ are given by (\ref{rec:eq:kernel-upper-value}%
)--(\ref{rec:eq:kernel-lower-value}).

\begin{proposition}
Under the hypotheses of Theorem~\ref{thm:transmutation-representation}, for $%
n\geq0$ and every $x\in(0,b]$,
\begin{equation}
\frac{2\chi_n^\pm(x)}{2n+1} =\sum_{k=0}^n l_{k,n}\frac{(\Phi_k^\pm)^{%
\prime}(x)}{x^k} -\frac{e^{\mp i\sigma(x)}}2 \left(\frac{n(n+1)}x+\mathcal{R}(x)\mp
ip(x)\right) \mp\frac{i(-1)^np(0)}2e^{\pm i\sigma(x)}.
\label{der:eq:direct-chi}
\end{equation}
\end{proposition}

\noindent\textit{Proof.} Differentiating (\ref{rec:eq:Phi-definition}) using
(\ref{rec:eq:transmutation-operator}), and substituting the characteristic
values (\ref{rec:eq:kernel-upper-value})--(\ref{rec:eq:kernel-lower-value}),
gives
\begin{equation}
(\Phi_k^\pm)^{\prime}(x) =e^{\mp i\sigma(x)} \left(kx^{k-1}+\frac{x^k}{2}(%
\mathcal{R}(x)\mp ip(x))\right) \pm\frac{ip(0)}2e^{\pm i\sigma(x)}(-x)^k +\int_{-x}^x K_{\pm,x}(x,t)t^k\,dt.
\label{der:eq:Phi-derivative}
\end{equation}%
Multiply (\ref{der:eq:Phi-derivative}) by $l_{k,n}x^{-k}$ and sum over $%
k=0,\ldots,n$. The identities
\begin{equation*}
\sum_{k=0}^n l_{k,n}=P_n(1)=1, \qquad \sum_{k=0}^n(-1)^k l_{k,n}=P_n(-1)=(-1)^n, \qquad \sum_{k=0}^n k l_{k,n}=P_n^{\prime}(1)=\frac{n(n+1)}2
\end{equation*}%
reduce the integral terms to $\int_{-x}^x K_{\pm,x}(x,t)P_n(t/x)\,dt$.
 \hfill$%
\square$

Alongside $\widetilde{\kappa}_n^\pm$ from (\ref%
{rec:eq:scaled-kappa-definition}), define the scaled derivative coefficients
\begin{equation*}
\widetilde{\chi}_n^\pm(x)=x^n\chi_n^\pm(x).
\end{equation*}

\begin{proposition}
The initial scaled derivative coefficients have the form
\begin{align}
\widetilde{\chi}_0^\pm(x) &=\frac{f^{\prime}(x)}2 -\frac{e^{\mp i\sigma(x)}}4(%
\mathcal{R}(x)\mp ip(x)) \mp\frac{ip(0)}4e^{\pm i\sigma(x)},
\label{der:eq:scaled-chi0}\\
\frac{2\widetilde{\chi}_1^\pm(x)}3 &=g^{\prime}(x)\mp ic_1^{\prime}(x) -e^{\mp i\sigma(x)} \left(1+\frac{x}{2}(\mathcal{R}(x)\mp ip(x))\right) \pm \frac{ixp(0)}2e^{\pm i\sigma(x)}.
\label{der:eq:scaled-chi1}
\end{align}%
Both functions vanish at $x=0$.
\end{proposition}

\noindent\textit{Proof.} For $n=0$, use $P_0=1$ and $\Phi_0^\pm=f$ in (\ref%
{der:eq:direct-chi}) to obtain (\ref{der:eq:scaled-chi0}). For $n=1$, $%
P_1(z)=z$, so (\ref{der:eq:direct-chi}) takes the form
\begin{equation*}
\begin{split}
\chi_1^\pm(x) =\frac32\left[ \frac{(\Phi_1^\pm)^{\prime}(x)}x -\frac{e^{\mp
i\sigma(x)}}2 \left(\frac2x+\mathcal{R}(x)\mp ip(x)\right) \pm\frac{ip(0)}%
2e^{\pm i\sigma(x)}\right].
\end{split}%
\end{equation*}%
Multiplication by $x$ and the identity $\Phi_1^\pm=g\mp ic_1$ give (\ref%
{der:eq:scaled-chi1}). \hfill$%
\square$

\begin{proposition}
\label{der:prop:initial-coefficients} The first derivative coefficients have
the form
\begin{align}
\widehat a_0(x)&=f^{\prime}(x) -\frac{\mathcal{R}(x)}2\cos\sigma(x) +\frac{%
p(x)+p(0)}2\sin\sigma(x),  \label{der:eq:ahat0}\\
\widehat u_0(x)&=\frac{\mathcal{R}(x)}2\sin\sigma(x) +\frac{p(x)-p(0)}%
2\cos\sigma(x),  \label{der:eq:uhat0}\\
\widehat v_0(x)&=3\left[ \frac{g^{\prime}(x)-\cos\sigma(x)}x +\frac{p(x)-p(0)}%
2\sin\sigma(x) -\frac{\mathcal{R}(x)}2\cos\sigma(x)\right],
\label{der:eq:vhat0}\\
\widehat b_0(x)&=3\left[ \frac{c_1^{\prime}(x)-\sin\sigma(x)}x -\frac{%
p(x)+p(0)}2\cos\sigma(x) -\frac{\mathcal{R}(x)}2\sin\sigma(x)\right].
\label{der:eq:bhat0}
\end{align}
\end{proposition}

\noindent\textit{Proof.} Substitute the initial coefficient formulas (\ref%
{der:eq:scaled-chi0})--(\ref{der:eq:scaled-chi1}) into the combinations (\ref%
{der:eq:ahat-bhat})--(\ref{der:eq:uhat-vhat}), using $\widetilde{\chi}%
_0^\pm=\chi_0^\pm$ and $\widetilde{\chi}_1^\pm=x\chi_1^\pm$. This gives (\ref%
{der:eq:ahat0})--(\ref{der:eq:bhat0}). The identities also hold at zero, since $\widehat v_0$ and $\widehat b_0$ are continuous there and the remaining terms have the indicated limits. \hfill$\square$

Thus, at $\rho=0$, formula (\ref{der:eq:c-derivative-nsbf}) gives $%
c^{\prime}(0,x)=f^{\prime}(x)$ by (\ref{der:eq:ahat0}), while (\ref%
{der:eq:s-derivative-nsbf}) gives $s^{\prime}(0,x)=g^{\prime}(x)$ by (\ref%
{der:eq:uhat0}) and (\ref{der:eq:vhat0}).

\subsection{Recurrent procedure for the derivative NSBF coefficients}

Recall $\mathcal{F}_\pm$ from (\ref{rec:eq:F-Legendre}) and, for fixed $x>0$%
, set
\begin{equation*}
\mathcal{H}_\pm(x,u)=xK_{\pm,x}(x,xu).
\end{equation*}
The derivative in $K_{\pm,x}$ is taken at fixed $t$, whereas $\mathcal{F}%
_{\pm,x}$ is taken at fixed $u$. Hence
\begin{equation}
\mathcal{H}_\pm =\mathcal{F}_{\pm,x} -\frac{\mathcal{F}_\pm+u\mathcal{F}%
_{\pm,u}}{x}.  \label{der:eq:H-F-chain-rule}
\end{equation}

\begin{proposition}
\label{prop:kernel-derivative-relation} Under the hypotheses of Theorem~\ref%
{thm:transmutation-representation}, the following identities hold for $n\ge2$
and every $x\in(0,b]$:
\begin{equation}
\chi_n^\pm =(\kappa_n^\pm)^{\prime }+\frac{n}{x}\kappa_n^\pm
-\frac{2n+1}{2n-3} \left[(\kappa_{n-2}^\pm)^{\prime }-\frac{n-1}{x}%
\kappa_{n-2}^\pm -\chi_{n-2}^\pm\right].
\label{der:eq:chi-local-relation}
\end{equation}
Equivalently, for the scaled coefficients,
\begin{equation}
\widetilde{\chi}_n^\pm ={}(\widetilde{\kappa}_n^\pm)^{\prime }-\frac{2n+1}{%
2n-3}x^2 \left[(\widetilde{\kappa}_{n-2}^\pm)^{\prime }-\frac{2n-3}{x}%
\widetilde{\kappa}_{n-2}^\pm -\widetilde{\chi}_{n-2}^\pm\right].
\label{der:eq:scaled-chi-local-differential}
\end{equation}
\end{proposition}

\noindent\textit{Proof.} Lemma~\ref{lem:goursat-sobolev-regularity} with $%
m=0 $ gives $\mathcal{F}_\pm(x,\cdot)\in W^{1,2}(-1,1)$ for every $x>0$. On
compact subintervals of $(0,b]$, the change of variables $t=xu$ preserves
the continuous $L^2$ dependence of $K_\pm$, $K_{\pm,x}$, and $K_{\pm,t}$
supplied by that lemma. Differentiating $\mathcal{F}_\pm(x,u)=xK_\pm(x,xu)$
with respect to $x$, with $u$ fixed, gives, for almost every $x>0$,
\begin{equation*}
\mathcal{F}_{\pm,x}(x,u) =K_\pm(x,xu)+xK_{\pm,x}(x,xu)+xuK_{\pm,t}(x,xu).
\end{equation*}%
The right-hand side is continuous with values in $L^2(-1,1)$. Consequently, $%
\mathcal{F}_\pm$ is continuously differentiable in $x$ with values in that
space, and the equality holds for every $x>0$ as an equality in $L^2(-1,1)$.
Its Legendre coefficients can therefore be differentiated under the integral
sign for every $x>0$.

For an index $m\geq 0$, integration by parts gives
\begin{equation*}
\begin{split}
\int_{-1}^{1}(\mathcal{F}_{\pm }+u\mathcal{F}_{\pm ,u})P_{m}\,du & =\left[ u\mathcal{F}_{\pm }(x,u)P_{m}(u)\right] _{-1}^{1}-%
\int_{-1}^{1}\mathcal{F}_{\pm }(x,u)uP_{m}^{\prime }(u)\,du \\
& =\mathcal{F}_{\pm }(x,1)+(-1)^{m}\mathcal{F}_{\pm
}(x,-1)-\int_{-1}^{1}\mathcal{F}_{\pm }(x,u)uP_{m}^{\prime }(u)\,du,
\end{split}%
\end{equation*}%
where we used $P_{m}(1)=1$ and $P_{m}(-1)=(-1)^{m}$. Taking a scalar product of both sides of (\ref%
{der:eq:H-F-chain-rule}) with $P_{m}$ and using (\ref{rec:eq:F-Legendre})
and (\ref{der:eq:chi-definition}) now yields
\begin{equation*}
\frac{2\left( \chi _{m}^{\pm }-(\kappa _{m}^{\pm })^{\prime }\right) }{2m+1}%
= -\frac{\mathcal{F}_{\pm }(x,1)+(-1)^{m}\mathcal{F}_{\pm }(x,-1)}{x}
 +\frac{1}{x}\int_{-1}^{1}\mathcal{F}_{\pm }(x,u)uP_{m}^{\prime }(u)\,du.
\end{equation*}%
Write this equality for $m=n$ and $m=n-2$ and subtract the second from the
first. The endpoint expressions are identical, so they cancel. By (\ref%
{rec:eq:u-Legendre-derivative-expansion}),
\begin{equation*}
u(P_{n}^{\prime }-P_{n-2}^{\prime })=nP_{n}+(n-1)P_{n-2}.
\end{equation*}%
Since $\int_{-1}^{1}\mathcal{F}_{\pm }P_{j}\,du=2\kappa _{j}^{\pm }/(2j+1)$,
the subtraction gives explicitly
\begin{equation*}
\begin{split}
\frac{2\left( \chi _{n}^{\pm }-(\kappa _{n}^{\pm })^{\prime }\right) }{2n+1%
}-\frac{2\left( \chi _{n-2}^{\pm }-(\kappa _{n-2}^{\pm })^{\prime }\right) }{%
2n-3} & =\frac{1}{x}\int_{-1}^{1}\mathcal{F}_{\pm }(x,u)u(P_{n}^{\prime
}(u)-P_{n-2}^{\prime }(u))\,du \\
&  =\frac{2n\kappa _{n}^{\pm }}{x(2n+1)}+\frac{2(n-1)\kappa _{n-2}^{\pm
}}{x(2n-3)}.
\end{split}%
\end{equation*}%
Multiplying by $(2n+1)/2$ and collecting terms proves (\ref%
{der:eq:chi-local-relation}). Finally, substitute
\begin{equation*}
(\kappa _{j}^{\pm })^{\prime }=x^{-j}\left[ (\widetilde{\kappa }_{j}^{\pm
})^{\prime }-\frac{j}{x}\widetilde{\kappa }_{j}^{\pm }\right] ,\qquad \chi
_{j}^{\pm }=x^{-j}\widetilde{\chi }_{j}^{\pm }
\end{equation*}%
for $j=n,n-2$ and multiply by $x^{n}$. This gives (\ref%
{der:eq:scaled-chi-local-differential}). \hfill $\square $

We now combine the local relation (\ref{der:eq:scaled-chi-local-differential}%
) with the solution recurrence (\ref{rec:eq:main-recursive-integration}).
The result contains no derivative of a previously computed coefficient.

\begin{theorem}
Assume $q\in L^2(0,b)$ and $p\in W^{1,2}(0,b)$, and let $F$ be a solution of
(\ref{rec:eq:auxiliary-equation}) without zeros on $[0,b]$, normalized by $%
F(0)=1$. For every $n\geq2$ and $0<x\leq b$,
\begin{equation}
\begin{split}
\frac{\widetilde\chi_n^\pm(x)}{2n+1} &=\frac{x^2\widetilde\chi_{n-2}^\pm(x)
-(2n-1)x\widetilde\kappa_{n-2}^\pm(x)}{2n-3} \\
&\quad+\frac{2(2n-1)}{2n-3} \left(F^{\prime }(x)\theta_n^\pm(x)+\frac{%
\eta_n^\pm(x)}{F(x)}\right) \mp2i \left(F^{\prime }(x)\mu_n^\pm(x)+\frac{\zeta_n^\pm(x)}{F(x)}%
\right).
\end{split}
\label{der:eq:scaled-chi-recursive-integration}
\end{equation}
The functions $\eta_{n}^\pm$, $\theta_{n}^\pm$, $\zeta_{n}^\pm$, and $%
\mu_{n}^\pm$ are defined by (\ref{rec:eq:eta-n})--(\ref{rec:eq:mu-n}).
\end{theorem}

\noindent \textit{Proof.} Set $C_{n}=(2n+1)/(2n-3)$ and $w=\widetilde{\kappa
}_{n-2}^{\pm }$. Differentiating (\ref{rec:eq:main-recursive-integration})
gives
\begin{equation*}
(\widetilde{\kappa }_{n}^{\pm })^{\prime }= C_{n}\{2xw+x^{2}w^{\prime
}+2(2n-1)[F^{\prime }\theta _{n}^{\pm }+F(\theta _{n}^{\pm })^{\prime }]\}  \mp 2i(2n+1)[F^{\prime }\mu _{n}^{\pm }+F(\mu _{n}^{\pm })^{\prime }].
\end{equation*}%
Substitute this equality into (\ref{der:eq:scaled-chi-local-differential}).
The two terms containing $x^{2}w^{\prime }$ cancel, and the remaining
expression has the form
\begin{equation*}
\widetilde{\chi }_{n}^{\pm }= C_{n}\{x^{2}\widetilde{\chi }_{n-2}^{\pm
}+(2n-1)xw+2(2n-1)[F^{\prime }\theta _{n}^{\pm }+F(\theta _{n}^{\pm
})^{\prime }]\}
\mp 2i(2n+1)[F^{\prime }\mu _{n}^{\pm }+F(\mu _{n}^{\pm })^{\prime }].
\end{equation*}%
By (\ref{rec:eq:theta-prime}) and (\ref{rec:eq:mu-prime}),
\begin{equation*}
F(\theta _{n}^{\pm })^{\prime }=\eta _{n}^{\pm }/F-x\widetilde{\kappa }%
_{n-2}^{\pm },\qquad F(\mu _{n}^{\pm })^{\prime }=\zeta _{n}^{\pm }/F.
\end{equation*}%
Substitution of these two identities gives (\ref%
{der:eq:scaled-chi-recursive-integration}). \hfill $\square $

Thus the same four successive integrals used for $\widetilde{\kappa}_n^\pm$
produce $\widetilde{\chi}_n^\pm$ by algebraic operations at each grid point,
with no additional quadrature or differentiation.

\begin{remark}The recurrence for the derivative coefficients can also be started at $n=1$ by introducing the auxiliary function
\[
\chi_{-1}^\pm(x)=-\chi_0^\pm(x)-K_\pm(x,-x)=-\chi_0^\pm(x)\mp\frac{i}{2}p(0)e^{\pm i\sigma(x)}.
\]
Together with $\kappa_{-1}^\pm=-\kappa_0^\pm$, this definition extends (\ref{der:eq:chi-local-relation}) and (\ref{der:eq:scaled-chi-recursive-integration}) to $n=1$, as follows directly from the initial coefficient formulas. Although $\widetilde{\chi}_{-1}^\pm=\chi_{-1}^\pm/x$ need not be bounded at zero, the integration formula uses it only through
\[
x^2\widetilde{\chi}_{-1}^\pm(x)=-x\bigl(\chi_0^\pm(x)+K_\pm(x,-x)\bigr),
\]
which extends continuously to zero with value zero. Thus $\chi_1^\pm$ requires no separate initialization.
\end{remark}

\subsection{Real-valued coefficients}

\label{der:sec:real-valued}

If $q$ and $p$ are real-valued, uniqueness for the conjugate Goursat
problems gives
\begin{equation}
K_- =\overline{K_+},\qquad \widetilde{\kappa}_n^-=\overline{\widetilde{\kappa%
}_n^+},\qquad \widetilde{\chi}_n^-=\overline{\widetilde{\chi}_n^+}.
\label{der:eq:real-conjugacy}
\end{equation}
Thus only the upper-sign recurrence needs to be computed. The coefficient
families can be recovered by
\begin{equation*}
\begin{split}
(a_m,u_m) & =2x^{-2m} (\operatorname{Re}\widetilde\kappa_{2m}^+, \operatorname{Im}%
\widetilde\kappa_{2m}^+), \\
(v_m,b_m) & =2x^{-2m-1} (\operatorname{Re}\widetilde\kappa_{2m+1}^+, -\operatorname{Im}%
\widetilde\kappa_{2m+1}^+).%
\end{split}%
\end{equation*}
The same relations recover the derivative coefficients: replace
$a_m,u_m,v_m,b_m$ by $\widehat a_m,\widehat u_m,\widehat v_m,\widehat b_m$,
respectively, and $\widetilde\kappa$ by $\widetilde\chi$.

\section{Computational procedure}
\label{sec:computation}

The numerical implementation follows the approach described in
\cite{KNT2017}. The computations in Subsection~\ref{num:sec:matlab-sech}
were performed in MATLAB R2017a using double precision. For the smooth
real coefficients in that example, the procedure is as follows.

\begin{enumerate}
\item On a uniform mesh of $[0,b]$, compute $f,g$ and their first
derivatives by the spectral parameter power series (SPPS) method
\cite{KravchenkoPorter2010}, and set $F=f+ig$. Evaluate $\sigma$ and
$\mathcal R$, as well as all subsequent integrals, by the six-point
Newton--Cotes rule used in \cite{KNT2017}.

\item After initialization, compute the scaled solution and derivative
coefficients together
using (\ref{rec:eq:main-recursive-integration}) and
(\ref{der:eq:scaled-chi-recursive-integration}). By linearity, the four
integrals in (\ref{rec:eq:eta-n})--(\ref{rec:eq:mu-n}) are combined into
two successive quadratures at each index. The derivative coefficients
use the same integrals. Only the upper-sign coefficients are computed; the lower
sign follows from (\ref{der:eq:real-conjugacy}).

\item Recover the unscaled coefficients at positive mesh nodes and
evaluate the finite spherical Bessel sums, retaining indices through
$N$. At $\rho=0$, use the corresponding limits, and at $x=0$ use the
prescribed initial conditions. The Dirichlet eigenvalues are found as
zeros of $s_N(\rho,b)$ by the MATLAB function
\texttt{fzero}, with \texttt{TolX} set to \texttt{eps}.
\end{enumerate}

The coefficients are computed once and reused for all spectral
parameters.

\section{Separated boundary conditions and spectral approximation}

\label{sec:spectral-approximation}

Consider the general separated conditions
\begin{equation}
A_{0}y(0)+A_{1}y^{\prime }(0)=0,\qquad B_{0}y(b)+B_{1}y^{\prime }(b)=0,
\label{spec:eq:general-bc}
\end{equation}%
where $A_{0},A_{1},B_{0},B_{1}$ are complex constants, with $%
(A_{0},A_{1})\neq (0,0)$ and $(B_{0},B_{1})\neq (0,0)$. The solution $%
A_{1}c(\rho ,x)-A_{0}s(\rho ,x)$ satisfies the left condition in (\ref%
{spec:eq:general-bc}). Hence a characteristic function of the problem (\ref%
{eq:pencil}), (\ref{spec:eq:general-bc}) can be written as follows
\begin{equation}
\Delta (\rho )=B_{0}\left( A_{1}c(\rho ,b)-A_{0}s(\rho ,b)\right)
+B_{1}\left( A_{1}c^{\prime }(\rho ,b)-A_{0}s^{\prime }(\rho ,b)\right) .
\label{spec:eq:characteristic-function}
\end{equation}%
The zeros of $\Delta $ are the values of the spectral parameter for which
the problem admits a nonzero solution. Dirichlet conditions correspond, for
example, to $(A_{0},A_{1})=(1,0)$ and $(B_{0},B_{1})=(1,0)$. For the Robin
conditions
\begin{equation*}
y^{\prime }(0)-h_{0}y(0)=0,\qquad y^{\prime }(b)+h_{b}y(b)=0,
\end{equation*}%
one takes $(A_{0},A_{1})=(-h_{0},1)$ and $(B_{0},B_{1})=(h_{b},1)$,
obtaining
\begin{equation*}
\Delta =c^{\prime }(\rho ,b)+h_{0}s^{\prime }(\rho ,b)+h_{b}\left( c(\rho
,b)+h_{0}s(\rho ,b)\right) .
\end{equation*}

We retain the terms with indices $0,\ldots ,N$ in the representations of the
two auxiliary solutions and their derivatives, and replace $c,s,c^{\prime
},s^{\prime }$ in (\ref{spec:eq:characteristic-function}) by $%
c_{N},s_{N},c_{N}^{[1]},s_{N}^{[1]}$ to define $\Delta _{N}$. Thus
\begin{equation}
\Delta _{N}(\rho )=B_{0}\left( A_{1}c_{N}(\rho ,b)-A_{0}s_{N}(\rho
,b)\right) +B_{1}\left( A_{1}c_{N}^{[1]}(\rho ,b)-A_{0}s_{N}^{[1]}(\rho
,b)\right) .  \label{spec:eq:truncated-characteristic-function}
\end{equation}%
To describe the error in the endpoint values, set
\begin{equation*}
E_{N}(b)=d_{+,N}(b)+d_{-,N}(b),\qquad \widehat{E}_{N}(b)=\widehat{d}%
_{+,N}(b)+\widehat{d}_{-,N}(b).
\end{equation*}

\begin{theorem}[Quantitative spectral approximation]
\label{thm:spectral-approximation} Assume
\begin{equation*}
q\in L^{2}(0,b),\qquad p\in W^{1,2}(0,b).
\end{equation*}%
Let $C\geq 0$ and $|\operatorname{Im}\rho |\leq C$. Then
\begin{equation}
|\Delta (\rho )-\Delta _{N}(\rho )|\leq \frac{|A_{1}|\Lambda
_{C}(b)+|A_{0}|\Theta _{C}(b)}{2}\left( |B_{0}|E_{N}(b)+|B_{1}|\widehat{E}%
_{N}(b)\right) =:\mathcal{E}_{N}(C).  \label{spec:eq:Delta-error}
\end{equation}%
Consequently, $\Delta _{N}\rightarrow \Delta $ uniformly on compact subsets
of the complex plane. If, moreover, $q\in W^{r,2}(0,b)$ and $p\in
W^{r+1,2}(0,b)$ for an integer $r\geq 1$, then
\begin{equation}
E_{N}(b)=O\left( N^{-r-1}\right) ,\qquad \widehat{E}_{N}(b)=O\left(
N^{-r}\right) ,  \label{spec:eq:coefficient-rates}
\end{equation}%
and hence $\Delta -\Delta _{N}=O(N^{-r})$ uniformly on each horizontal
strip. When $B_{1}=0$, this improves to $O(N^{-r-1})$.

Let $\Gamma$ be a positively oriented simple closed rectifiable contour
contained in the strip. If $\Delta$ has no zero on $\Gamma$ and either
\begin{equation*}
\mathcal{E}_N(C)<\min_{\rho\in\Gamma}|\Delta(\rho)|
\end{equation*}
or an available upper bound $\widetilde{\mathcal{E}}_N(C)\geq\mathcal{E}%
_N(C) $ satisfies
\begin{equation*}
\widetilde{\mathcal{E}}_N(C)< \min_{\rho\in\Gamma}|\Delta_N(\rho)|,
\end{equation*}
then $\Delta$ and $\Delta_N$ have the same number of zeros inside $\Gamma$,
counted with multiplicity. In particular, if $\rho_0$ is a simple eigenvalue
and a closed disk centered at $\rho_0$ contains no other eigenvalue, then
that disk contains exactly one zero of $\Delta_N$ for all sufficiently large
$N$.
\end{theorem}

\noindent\textit{Proof.} Define
\begin{equation*}
\delta c=c-c_N,\quad \delta s=s-s_N,\quad \delta c^{\prime }=c^{\prime
}-c_N^{[1]},\quad \delta s^{\prime }=s^{\prime }-s_N^{[1]}.
\end{equation*}
Subtracting (\ref{spec:eq:truncated-characteristic-function}) from (\ref%
{spec:eq:characteristic-function}) and using the triangle inequality, we
obtain
\begin{equation*}
|\Delta-\Delta_N| \leq |B_0|\left(|A_1||\delta c|+|A_0||\delta s|\right) +|B_1|\left(|A_1||\delta c^{\prime }|+|A_0||\delta s^{\prime }|\right).
\end{equation*}
The estimates (\ref{err:eq:c-estimate}), (\ref{err:eq:s-estimate}), (\ref%
{der:err:c-prime}), and (\ref{der:err:s-prime}) now give (\ref%
{spec:eq:Delta-error}). Every compact subset of the $\rho$-plane is
contained in a horizontal strip. At $x=b$, completeness of the Legendre
polynomials gives $E_N(b)\to0$ and $\widehat E_N(b)\to0$. It follows from (%
\ref{spec:eq:Delta-error}) that $\Delta_N\to\Delta$ uniformly on compact
subsets.

Suppose in addition that $q\in W^{r,2}(0,b)$ and $p\in W^{r+1,2}(0,b)$.
Theorem~\ref{thm:kernel-projection-rates} may then be applied with $r+1$,
because its hypotheses become precisely $q\in W^{r,2}(0,b)$ and $p\in
W^{r+1,2}(0,b)$. It gives $E_N(b)=O(N^{-r-1})$. Theorem~\ref%
{thm:derivative-error-estimates}, applied with $r$, gives $\widehat
E_N(b)=O(N^{-r})$. Substitution into (\ref{spec:eq:Delta-error}) proves (\ref%
{spec:eq:coefficient-rates}), the $O(N^{-r})$ estimate for general separated
conditions, and the $O(N^{-r-1})$ estimate when $B_1=0$.

Under the first contour criterion, (\ref{spec:eq:Delta-error}) gives $%
|\Delta_N-\Delta|<|\Delta|$ on $\Gamma$. Rouch\'e's theorem, with $\Delta$
as the comparison function, shows that $\Delta_N$ and $\Delta$ have the same
number of zeros inside $\Gamma$. The second criterion requires an upper bound
for the truncation error, which may be obtained, for example, from quantitative
regularity estimates for the kernels.  Under the second criterion, $|\Delta_N-\Delta|\leq%
\mathcal{E}_N(C)\leq \widetilde{\mathcal{E}}_N(C)<|\Delta_N|$ on $\Gamma$,
and Rouch\'e's theorem is applied with $\Delta_N$ as the comparison function.

Finally, let $\rho_0$ be a simple zero of $\Delta$. Choose a closed disk
centered at $\rho_0$ whose boundary contains no zero of $\Delta$ and whose
interior contains no other zero. Uniform convergence on the boundary and the
first contour criterion imply that, for all sufficiently large $N$, the disk
contains exactly one zero of $\Delta_N$, counted with multiplicity. The disk
may be chosen with arbitrarily small radius. Hence these zeros converge to $%
\rho_0$. \hfill$\square$

\begin{remark}
For the Robin conditions, (\ref{spec:eq:Delta-error}) becomes
\begin{equation*}
|\Delta -\Delta _{N}|\leq \frac{\Lambda _{C}(b)+|h_{0}|\Theta _{C}(b)}{2}%
\left( |h_{b}|E_{N}(b)+\widehat{E}_{N}(b)\right) .
\end{equation*}%
On the real spectral axis the right-hand side is independent of $|\rho |$.
Thus the estimate does not deteriorate as the real spectral parameter
increases.
\end{remark}

\section{Numerical examples}
\label{sec:example}

\subsection{An exactly solvable example}

This example tests the solution and eigenvalue approximations against closed
formulas. Because $p\ne0$, it retains the spectral dependence of the
quadratic pencil. We compare the solutions throughout a broad real spectral
interval and check the Dirichlet eigenvalues on both branches.

Let $b=1$, $q(x)=0$, and $p(x)=-1/2$. The fundamental solutions of (\ref%
{eq:pencil}) have the form
\begin{equation*}
c(\rho ,x)=\cos (\nu x),\qquad s(\rho ,x)=\frac{\sin (\nu x)}{\nu },\qquad
\nu ^{2}=\rho ^{2}+\rho .
\end{equation*}%
The explicit solutions also provide derivatives for checking their NSBF
series. Below we compare the solutions and eigenvalues. For Dirichlet
conditions at $0$ and $1$, both branches of the spectrum can be written
explicitly:
\begin{equation}
\rho _{n}^{\pm }=\frac{-1\pm \sqrt{1+4n^{2}\pi ^{2}}}{2},\qquad n=1,2,\ldots
.  \label{num:eq:exact-spectrum}
\end{equation}%
For instance, (\ref{num:eq:exact-spectrum}) gives
\begin{equation*}
\begin{array}{c|rr}
n & \rho _{n}^{+} & \rho _{n}^{-} \\ \hline
1 & 2.6811325658 & -3.6811325658 \\
2 & 5.8030482788 & -6.8030482788 \\
3 & 8.9380315538 & -9.9380315538 \\
40 & 125.1647008580 & -126.1647008580%
\end{array}%
\end{equation*}%
To check the series, use (\ref{rec:eq:direct-kappa-formula}) at $x=1$ and
truncate both sign series after index $N$. The moments needed in that
formula follow from $\Phi _{k}^{+}(1)=k![z^{k}](\cos \nu +z\sin \nu /\nu )$,
where $[z^{k}]F(z)$ denotes the coefficient of $z^{k}$ in the power series
expansion of $F(z)$, $z=i\rho $, and $\nu ^{2}=-z^{2}-iz$. For real $\rho $
we use complex conjugation for the lower sign. On 800 equally spaced real
values of $\rho $ in $[1,200]$, we measure the solution error by
\begin{equation*}
E_{N}=\max_{\rho }\max \{|c(\rho ,1)-c_{N}(\rho ,1)|,|\rho (s(\rho
,1)-s_{N}(\rho ,1))|\}.
\end{equation*}%
For each of the first 40 eigenvalues on either branch, locate the zero of $%
s_{N}(\rho ,1)$ nearest the exact value and denote the largest eigenvalue
error by $D_{N}$. Evaluation of the alternating direct coefficient formula
with 50-digit arithmetic gives the values in the table
\begin{equation*}
\begin{array}{c|rr}
N & E_{N} & D_{N} \\ \hline
2 & 4.02\cdot 10^{-3} & 3.41\cdot 10^{-3} \\
4 & 3.89\cdot 10^{-5} & 4.06\cdot 10^{-5} \\
6 & 1.90\cdot 10^{-7} & 2.01\cdot 10^{-7} \\
8 & 5.65\cdot 10^{-10} & 4.98\cdot 10^{-10}%
\end{array}%
\end{equation*}%
Both branches are unbounded. The solutions are entire functions of $\nu
^{2}=\rho ^{2}+\rho $, so their values are independent of the square-root
branch.

\subsection{A pencil with a hyperbolic-secant potential}
\label{num:sec:matlab-sech}

Consider the equation \cite{Crum1955,Lekner2007}
\[
 -y''+\left(-\frac{2}{\cosh^2(x-1)}+\rho\right)y=\rho^2y,
 \qquad 0\leq x\leq2,
\]
that is, $p(x)=1/2$ and $q(x)=-2/\cosh^2(x-1)$. In this example the
potential $q$ is nonconstant, and the fundamental solutions admit explicit
expressions with which the NSBF approximations can be compared.

Denote $\lambda=\rho^2-\rho$, $T(x)=\tanh(x-1)$, $a=\tanh1$, and
$r_0=1-a^2$. If $v''+\lambda v=0$, then $y=v'-Tv$ satisfies
$-y''-2y/\cosh^2(x-1)=\lambda y$. Choosing the linear combinations
with the prescribed initial conditions gives
\begin{equation}
\begin{aligned}
c(\rho,x)
&=\frac{(\lambda+r_0-aT(x))\cos(kx)
-\bigl(\lambda a+(\lambda+r_0)T(x)\bigr)\sin(kx)/k}{\lambda+1},\\
s(\rho,x)
&=\frac{(T(x)+a)\cos(kx)
+(\lambda-aT(x))\sin(kx)/k}{\lambda+1},
\qquad k^2=\lambda=\rho^2-\rho.
\end{aligned}
\label{num:eq:matlab-sech-solutions}
\end{equation}
Here $\sin(kx)/k=x$ when $k=0$. Both expressions are even in $k$,
and hence do not depend on the choice of the square root. The point
$\lambda=-1$ does not occur on the spectral sets considered below.

For the boundary conditions $y(0)=y(2)=0$, the characteristic function is
\[
 s(\rho,2)=\frac{2a\cos(2k)+(k^2-a^2)\sin(2k)/k}{k^2+1}.
\]
Let $\lambda_n$ be the Dirichlet eigenvalues of
$-y''-2y/\cosh^2(x-1)=\lambda y$ on $[0,2]$. Comparing $q$ with the constant potentials $-2$ and $0$ and applying \cite[Corollary~4.13]{Teschl2014} gives the following bounds
\[
 \left(\frac{n\pi}{2}\right)^2-2
 \leq\lambda_n\leq\left(\frac{n\pi}{2}\right)^2,
 \qquad n=1,2,\ldots.
\]
In particular, $\lambda_n>0$.

Writing the numerator of the characteristic
function as a single sine gives explicit expressions for the eigenvalues 
\begin{equation}
 k_n+\arctan\frac{a}{k_n}=\frac{n\pi}{2},\qquad
 \lambda_n=k_n^2,\qquad
 \rho_n^\pm=\frac12\pm\sqrt{\frac14+\lambda_n}.
\label{num:eq:matlab-sech-spectrum}
\end{equation}

The numerical computation was performed in MATLAB R2017a using double precision. We used $10001$ equally spaced points on
$[0,2]$. The functions $f,g$ and their first derivatives were computed using the
SPPS representation \cite{KravchenkoPorter2010}, with formal powers through index
$64$, and $F=f+ig$ was used as a non-vanishing particular solution in the integral recurrences. We have computed the solution and
derivative coefficients together through index $24$, using the six-point Newton--Cotes procedure described in
Section~\ref{sec:computation}.

For each $N$, denoting the largest retained Bessel index, the first $200$ eigenvalues
on each branch were found as zeros of $s_N(\rho,2)$ by MATLAB's function \texttt{fzero},
with \texttt{TolX} set to \texttt{eps}. The intervals used by \texttt{fzero} were
obtained from the bounds for $\lambda_n$ above. The eigenvalue errors were evaluated
using reference values computed with higher precision from (\ref{num:eq:matlab-sech-spectrum}).

We measured the solution and derivative errors on the real spectral grid
\[
 \mathcal G=\{-200+j/4:j=0,\ldots,1600\}\cup\{-10^{-8},10^{-8}\},
\]
by
\[
\begin{aligned}
 E_N&=\max_{\rho\in\mathcal G}\max\{|c_N(\rho,2)-c(\rho,2)|, |\rho(s_N(\rho,2)-s(\rho,2))|\},\\
 E_N^{(1)}&=\max_{\rho\in\mathcal G}\max\{|c_N^{[1]}(\rho,2)-c'(\rho,2)|,|s_N^{[1]}(\rho,2)-s'(\rho,2)|\}.
\end{aligned}
\]
Here $c_N^{[1]}$ and $s_N^{[1]}$ are the approximations obtained from
the derivative NSBF representations. We also evaluated the solution
error on the two finite complex sets
\[
 \{-200+j/2+i/4:j=0,\ldots,800\},\qquad \{-200+j/2-i/4:j=0,\ldots,800\}.
\]
Let $E_N^{\rm strip}$ be defined as $E_N$, with $\mathcal G$ replaced
by the union of these two sets. Finally, if $\rho_{n,N}^\pm$ are the
computed eigenvalues, define
\[
 D_N=\max_{1\leq n\leq200}\max\{|\rho_{n,N}^+-\rho_n^+|,|\rho_{n,N}^--\rho_n^-|\}.
\]

Table~\ref{num:tab:matlab-sech-errors} presents the errors obtained in
MATLAB. Increasing $N$ beyond $16$ gives little further improvement on
this mesh. The solution and eigenvalue errors are approximately
$3.3\cdot10^{-13}$ and $2.3\cdot10^{-13}$, respectively, while the
absolute derivative error is approximately $1.6\cdot10^{-11}$.

\begin{table}[ht]
\centering
\caption{Errors for the hyperbolic-secant potential, computed in MATLAB
R2017a using $10001$ mesh points.}
\label{num:tab:matlab-sech-errors}
\begin{tabular}{r|rrrr}
$N$ & $E_N$ & $E_N^{(1)}$ & $E_N^{\rm strip}$ & $D_N$\\
\hline
4  & $5.06\cdot10^{-3}$  & $1.75\cdot10^{-3}$
   & $5.29\cdot10^{-3}$  & $2.42\cdot10^{-3}$\\
8  & $1.34\cdot10^{-6}$  & $1.21\cdot10^{-6}$
   & $1.38\cdot10^{-6}$  & $7.39\cdot10^{-7}$\\
12 & $6.24\cdot10^{-11}$ & $2.44\cdot10^{-11}$
   & $6.29\cdot10^{-11}$ & $2.55\cdot10^{-11}$\\
16 & $3.28\cdot10^{-13}$ & $1.60\cdot10^{-11}$
   & $3.44\cdot10^{-13}$ & $2.26\cdot10^{-13}$\\
20 & $3.28\cdot10^{-13}$ & $1.60\cdot10^{-11}$
   & $3.44\cdot10^{-13}$ & $2.26\cdot10^{-13}$\\
24 & $3.28\cdot10^{-13}$ & $1.59\cdot10^{-11}$
   & $3.44\cdot10^{-13}$ & $2.26\cdot10^{-13}$\\
\hline
\end{tabular}
\end{table}

Some of the eigenvalues computed with $N=16$ are listed in
Table~\ref{num:tab:matlab-sech-eigenvalues}. The values are printed to
ten decimal places.

\begin{table}[ht]
\centering
\caption{Selected Dirichlet eigenvalues computed in MATLAB with $N=16$.}
\label{num:tab:matlab-sech-eigenvalues}
\begin{tabular}{r|rr}
$n$ & $\rho_{n,16}^+$ & $\rho_{n,16}^-$\\
\hline
1   & $1.4654239069$   & $-0.4654239069$\\
2   & $3.4263868257$   & $-2.4263868257$\\
3   & $5.0738262471$   & $-4.0738262471$\\
40  & $63.3217199683$  & $-62.3217199683$\\
200 & $314.6572390057$ & $-313.6572390057$\\
\hline
\end{tabular}
\end{table}

\begin{figure}[ht]
\centering
\includegraphics[width=.88\textwidth]{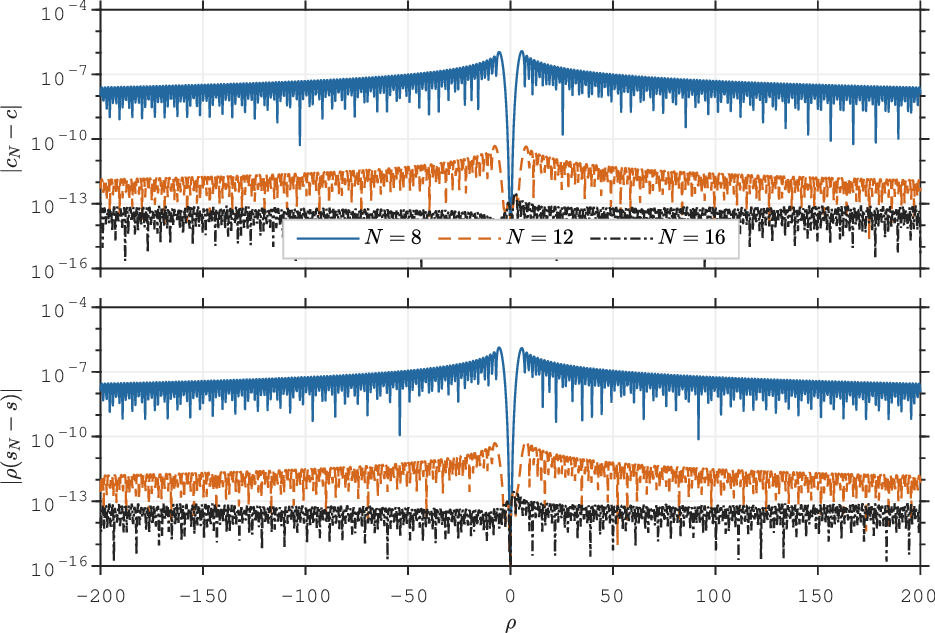}
\caption{Errors $|c_N(\rho,2)-c(\rho,2)|$ (upper panel) and
$|\rho(s_N(\rho,2)-s(\rho,2))|$ (lower panel) for $N=8,12,16$,
computed in MATLAB.}
\label{num:fig:matlab-sech-solutions}
\end{figure}

\begin{figure}[ht]
\centering
\includegraphics[width=.88\textwidth]{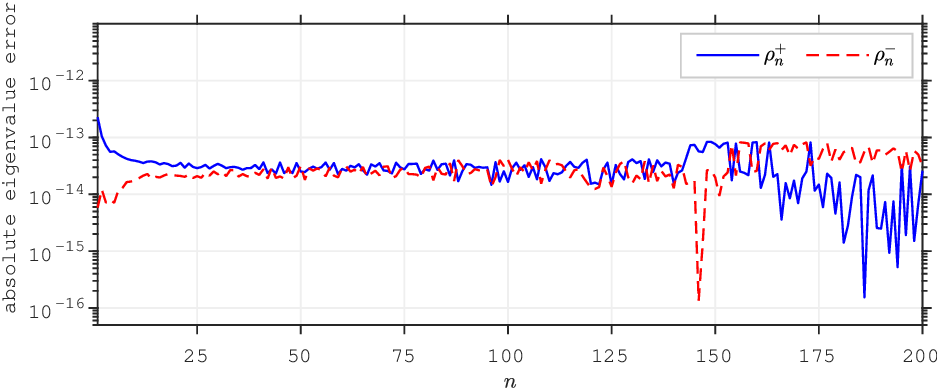}
\caption{Absolute errors in the first $200$ Dirichlet eigenvalues on
each branch for $N=16$, computed in MATLAB.}
\label{num:fig:matlab-sech-eigenvalues}
\end{figure}

Figure~\ref{num:fig:matlab-sech-solutions} shows the two solution errors
on the real spectral grid for $N=8,12,16$.
Figure~\ref{num:fig:matlab-sech-eigenvalues} shows the individual
eigenvalue errors on both branches for $N=16$. The fluctuations for
large indices are consistent with double-precision rounding.

\section{Conclusions}

In the present work, representations in the form of Neumann series of Bessel
functions are obtained for the solutions of the quadratic pencil and their
first derivatives. Their coefficients can be calculated by a simple procedure involving successive integrations of known functions. The truncation error estimates remain
uniform as the real part of the spectral parameter grows and give
corresponding error control when computing spectra of boundary value
problems with separated boundary conditions. Two examples illustrate
the accuracy of the solution and eigenvalue approximations on a wide
spectral range.

The first NSBF coefficients give local recovery formulas for both
potentials, with complementary conditions of validity. This feature of the
representations is valuable for solving the corresponding inverse problems,
a topic for future work.

\section*{Data availability}
The data that support the findings of this study are available upon reasonable request.

\section*{Conflict of interest}
This work does not have any conflict of interest.

\section*{Declaration of AI use}

During the preparation of this manuscript, the authors used ChatGPT
(OpenAI) to assist with literature searches, numerical examples and revision of the text. The authors take full responsibility for the
content of the manuscript.

\appendix

\section{Proof of Lemma 2.2}

\label{app:kernel-regularity}

\noindent\textit{Proof.} Write $r=m+1$. The one dimensional Sobolev
embedding and product rule imply that $W^{k,2}(0,b)$ is an algebra for each
integer $k\geq1$ \cite[Chapter~5]{AdamsFournier}. Applied to (\ref{eq:phase}%
), (\ref{eq:goursat-boundary-functions}), and the identity $\mathcal{R}%
^{\prime}=q+p^2$, this gives $h_{\pm,1},h_{\pm,2}\in W^{r,2}(0,b)$. Their
norms depend locally on $\|q\|_{W^{r-1,2}(0,b)}+\|p\|_{W^{r,2}(0,b)}$, and
the map from $(q,p)$ to these data is locally Lipschitz in the same norms.
For the analogous regularity of kernels with real coefficients, see \cite[%
Theorem~2.1]{Guseinov2014}.

First let us consider smooth coefficients and denote the right side of (\ref%
{eq:goursat-H}) by $G_{\pm }$. Integrating in each coordinate gives
\begin{equation}
\begin{split}
H_{\pm ,\xi }(\xi ,\eta )& =h_{\pm ,1}^{\prime }(\xi )+\int_{0}^{\eta
}G_{\pm }(\xi ,v)\,dv, \\
H_{\pm ,\eta }(\xi ,\eta )& =h_{\pm ,2}^{\prime }(\eta )+\int_{0}^{\xi
}G_{\pm }(u,\eta )\,du.
\end{split}
\label{eq:goursat-first-derivatives}
\end{equation}%
For the restrictions of the derivatives to the line $\xi+\eta=x$, write
\begin{equation*}
M_{k}(U;x)=\sum_{j+\ell \leq k}\left( \int_{0}^{x}\left\vert (\partial _{\xi
}^{j}\partial _{\eta }^{\ell }U)(s,x-s)\right\vert ^{2}\,ds\right)
^{1/2},\qquad 0<x\leq b.
\end{equation*}%
For $v\in L^{1}(0,b)$, Minkowski's inequality gives, for example,
\begin{equation*}
\left\Vert \int_{0}^{x-\xi }v(\xi +z)U(\xi ,z)\,dz\right\Vert
_{L^{2}(0,x;d\xi )}\leq \int_{0}^{x}|v(s)|\left(
\int_{0}^{s}|U(u,s-u)|^{2}\,du\right) ^{1/2}ds.
\end{equation*}%
The other first order integral in (\ref{eq:goursat-volterra}) has the same
estimate, and its double integral has this bound multiplied by $b$. Apply
these estimates first to (\ref{eq:goursat-volterra}), then to (\ref%
{eq:goursat-first-derivatives}). For higher $k$, differentiate the latter
identities $k-1$ times. The terms containing order $k$ derivatives of $%
H_{\pm }$ stay inside a Volterra integral. Mixed derivatives come directly
from (\ref{eq:goursat-H}). The product rule and the one-dimensional
embedding control the lower order terms, including the highest coefficient
derivatives in $L^{2}$. Induction gives the inequality
\begin{equation}
M_{k}(H_{\pm };x)\leq C_{k}+C_{k}\int_{0}^{x}(1+|p(s)|+|q(s)|)M_{k}(H_{\pm
};s)\,ds,\qquad 0\leq k\leq r.  \label{eq:goursat-induction-estimate}
\end{equation}%
Each $C_{k}$ is locally bounded in the stated coefficient norms and in the
bounds already obtained for $M_{0},\ldots ,M_{k-1}$. We use the following
integral form of Gronwall's inequality \cite[Chapter~III, Section~1]%
{HartmanODE}: if $h\geq 0$ belongs to $L^{1}(0,b)$ and $u(x)\leq
a+\int_{0}^{x}h(s)u(s)\,ds$, then $u(x)\leq a\exp (\int_{0}^{x}h(s)\,ds)$.
Apply it successively to (\ref{eq:goursat-induction-estimate}), taking $%
h(s)=C_{k}(1+|p(s)|+|q(s)|)$. This bounds $M_{r}(H_{\pm };x)$ uniformly for $%
0<x\leq b$.

Take smooth $q_{\nu }\rightarrow q$ in $W^{r-1,2}(0,b)$ and $p_{\nu
}\rightarrow p$ in $W^{r,2}(0,b)$. Apply the same induction to the
difference of their Goursat problems. The differing boundary terms
are controlled by the local Lipschitz estimate above. The resulting $M_{r}$
bound is at most a constant times $\Vert q_{\nu }-q_{\mu }\Vert
_{W^{r-1,2}}+\Vert p_{\nu }-p_{\mu }\Vert _{W^{r,2}}$, uniformly in $x$.
Thus every derivative along $\xi +\eta =x$, extended by zero to $(0,b)$,
converges uniformly in $x$ with values in $L^{2}(0,b)$. The change of
variables $x=\xi +\eta $ gives the equality
\begin{equation*}
\int_{\Omega _{b}}|\partial _{\xi }^{j}\partial _{\eta }^{\ell }U|^{2}\,d\xi
\,d\eta =\int_{0}^{b}\int_{0}^{x}|(\partial _{\xi }^{j}\partial _{\eta
}^{\ell }U)(s,x-s)|^{2}\,ds\,dx.
\end{equation*}%
Consequently the smooth kernels converge in $W^{r,2}(\Omega _{b})$. The
uniform convergence from (\ref{eq:goursat-volterra}) identifies this limit
with $H_{\pm }$ and preserves (\ref{eq:goursat-data}). The restrictions
of the smooth kernels and their derivatives to $\xi+\eta=x$, extended by zero, are continuous as $%
L^{2}$-valued functions of $x$ and converge uniformly in that space. Their
limits are continuous, and weak differentiation on $(-x,x)$ identifies the
limits as the stated $t$ derivatives. The identities $K_{\pm ,t}=(H_{\pm
,\xi }-H_{\pm ,\eta })/2$ and $K_{\pm ,x}=(H_{\pm ,\xi }+H_{\pm ,\eta })/2$,
together with $dt=2\,d\xi $, prove
(\ref{eq:goursat-section-bound})--(\ref{eq:goursat-section-bound-derivative}) and the
asserted continuity. \hfill $\square $

\end{document}